\documentclass{aims}
\usepackage{amsmath}
  \usepackage{paralist}
  \usepackage{graphics} %% add this and next lines if pictures should be in esp format
  \usepackage{epsfig} %For pictures: screened artwork should be set up with an 85 or 100 line screen
\usepackage{graphicx}  \usepackage{epstopdf}%This is to transfer .eps figure to .pdf figure; please compile your paper using PDFLeTex or PDFTeXify.
 \usepackage[colorlinks=true]{hyperref}
\hypersetup{urlcolor=blue, citecolor=red}

\def\currentvolume{X}
 \def\currentissue{X}
  \def\currentyear{202X}
   \def\currentmonth{XX}
    \def\ppages{X--XX}

\theoremstyle{definition}

\usepackage{multirow,xfrac}
\usepackage{graphics}
\usepackage{subcaption}
\usepackage[table]{xcolor}

\newcommand{\bv}{\pmb{v}}
\newcommand{\bw}{\pmb{w}}

\newcommand{\bn}{\pmb{n}}

\newcommand{\bx}{\pmb{x}}

\newcommand{\bD}{\pmb{D}}

\newcommand{\bu}{\pmb{u}}

\newcommand{\iK}{\mathcal{K}}

\newcommand{\mR}{\mathbb{R}}

\newcommand{\Div}{\text{div}\;}
\newcommand{\iS}{\mathcal{S}}
\newcommand{\iT}{\mathcal{T}}

\newcommand{\iP}{\mathcal{P}}

\newcommand{\iE}{\mathcal{E}}

\newcommand{\bphi}{\pmb{\varphi}}
\newcommand{\bphiG}{\pmb{\varphi}}
\newcommand{\bpsi}{\pmb{\phi}}

\title[Ventcel-Schwarz waveform relaxation in mixed formulations] %Use the shortened version of the full title
      {Optimized Ventcel-Schwarz waveform relaxation and mixed hybrid finite element method for transport problems}

\author[T-T.-P.  Hoang]{}


 \email{tzh0059@auburn.edu}

\thanks{Dedicated to Professor Georg Hetzer on the occasion of his $75^{\text{th}}$ birthday}
\thanks{This work is partially supported by the US National Science Foundation under grant number DMS-1912626.}

\begin{document}
\maketitle

% Enter the first author's name and address:
\centerline{\scshape Thi-Thao-Phuong Hoang$^*$}
\medskip
{\footnotesize
% please put the address of the first author
 \centerline{Department of Mathematics and Statistics, Auburn University}
   %\centerline{Other lines}
   \centerline{Auburn, AL 36849, USA}
} % Do not forget to end the {\footnotesize by the sign }

\medskip

\bigskip

% The name of the associate editor will be entered by an editorial staff
% "Communicated by the associate editor name" is not needed for special issue.
% \centerline{(Communicated by the associate editor name)}

%The abstract of your paper
\begin{abstract}
%This is the abstract of your paper and it should not exceed \textbf{200} words.
This paper is concerned with the optimized Schwarz waveform relaxation method and Ventcel transmission conditions for the linear advection-diffusion equation.  A mixed formulation is considered in which the flux variable represents both diffusive and advective flux, and Lagrange multipliers are introduced on the interfaces between nonoverlapping subdomains to handle tangential derivatives in the Ventcel conditions.  A space-time interface problem is formulated and is solved iteratively. Each iteration involves the solution of time-dependent problems with Ventcel boundary conditions in the subdomains. The subdomain problems are discretized in space by a mixed hybrid finite element method based on the lowest-order Raviart-Thomas space and in time by the backward Euler method. The proposed algorithm is fully implicit and enables different time steps in the subdomains.  Numerical results with discontinuous coefficients and various Pecl\'et numbers validate the accuracy of the method with nonconforming time grids and confirm the improved convergence properties of Ventcel conditions over Robin conditions. 

\end{abstract}

%
% ------------------------------------------------
%
%    SECTION 1 : Introduction
%	
% ------------------------------------------------
%
\section{Introduction}
\label{Sec:intro}  
%\begin{itemize}
%\item Optimized Schwarz and Optimized Schwarz Waveform Relaxation methods
%\item Mixed formulations: upwind-mixed hybrid finite element method
%\item Previous work: Robin transmission conditions
%$\longrightarrow$ This work: Ventcel condition for both steady and unsteady problems (algorithm development and numerical experiments)
%\end{itemize}

Optimized Schwarz waveform relaxation (OSWR) methods are a class of global-in-time domain decomposition methods for parallel solutions of evolution problems.  They are iterative algorithms that solve time-dependent problems in the subdomains over the whole time interval and exchange data on the space-time interfaces through transmission operators of Robin or Ventcel~\cite{Ventcel} types.  Ventcel transmission conditions are second-order differential conditions which involve time and tangential derivatives; these conditions were first introduced for steady convection-diffusion problems in \cite{NR95,Japhet98}.  Differently from the classical Schwarz waveform relaxation method which exchanges only Dirichlet data on the interfaces and requires the subdomains to overlap,  OSWR methods converge with or without overlap.  The Robin or Ventcel transmission conditions include some coefficients that are determined by optimization of the convergence factor, thus the convergence of OSWR is significantly more improved than the classical approach. 
%can be optimized to improve the convergence rates of the iterations.  
OSWR methods were first introduced for the advection-reaction-diffusion and wave equations in one dimension in \cite{GHN99,OSWRwave03},  then extended to two dimensional convection-diffusion problems in \cite{Martin05}.  Analysis of the optimization problems was carried out in \cite{GH07} and \cite{BGH09} for the Robin and Ventcel transmission conditions applied to the one dimensional advection-diffusion equations; the two-dimensional case was studied in \cite{BGGH16}.  In addition to enhanced convergence properties,  OSWR allows different discretizations in both space and time in the subdomains,  which makes the methods well-suited for heterogeneous and coupled problems.  In \cite{GHK07,BlayoHJ07,HJO10,Japhet12,BDF13,HHM13},  discontinuous coefficients and nonmatching time discretizations were considered where a suitable time projection was employed to exchange information between the subdomains on the space-time interfaces.  The method was also applied to the viscous primitive equations of the ocean in \cite{Merlet10}.  Reviews of OSWR methods can be found in \cite{HalpernDD17, Gander08}.

%These iterative methods involve the solution of time-dependent problems in the subdomains at each iteration and information exchange on space-time boundary data through general (Robin or Ventcel) transmission operators in which coefficients can be optimized to improve convergence rates.  
%They are a combination of the optimized Schwarz method for steady problems \cite{Japhet01,GMN02, Gander06} and the classical waveform relaxation algorithm \cite{WR82} for systems of ordinary differential equations.  

For flow and transport problems in porous media,  it is important to use conservative cell-centered techniques such as mixed methods \cite{BF91,RT91} to obtain accurate approximations of the solutions.  %with their mass conservation property and simultaneous approximation of the same order of both the scalar and vector fields. %a natural way to handle heterogeneous and anisotropic tensors.  %In addition,  these methods give an approximation simultaneously, and to the same order, of both the scalar and vector field.  
The OSWR methods have been extensively studied mostly for the primal formulation with either Lagrange finite element or finite volume discretizations.  In the context of mixed formulations, OSWR methods with Robin transmission conditions and nonconforming time grids were studied for pure diffusion problems in \cite{H13} and for the advection-diffusion problems in \cite{H17,H21}.  Operator splitting was used in \cite{H17} so that the advection is treated explicitly and the diffusion implicitly,  while the method in \cite{H21} is fully implicit in time and the problem is discretized in space by mixed hybrid finite elements~\cite{Radu11, Radu14}.  In \cite{H16},  OSWR methods were applied to a reduced fracture model of the flow of a compressible fluid in a porous medium in which the fracture is treated as an interface between two subdomains and the so-called Ventcel-to-Robin transmission conditions were derived for such a model.  In \cite{H16Ventcel},  optimized Schwarz methods with Ventcel conditions in mixed form were considered for the steady diffusion problems.  It should be noted that the convergence of the OSWR method with Ventcel transmission conditions is improved over that with Robin conditions as shown in~\cite{BGH09,BGGH16} (where the problem is written in primal form).  

The objective of this work is to develop a global-in-time optimized Ventcel-Schwarz method for mixed formulations of the advection-diffusion problem.  The method is based on OSWR with Ventcel conditions, though it is not obtained in such a straightforward manner as in the case of primal formulations.  In particular,  Lagrange multipliers have to be introduced on the interfaces to handle tangential derivatives involved in the Ventcel conditions.  We consider nonoverlapping subdomains and formulate the the initial-boundary value problem on the whole domain as a space–time interface problem, through the use of trace operators.  Such an interface problem is solved iteratively in which each iteration involves solution of time-dependent subdomain problems with Ventcel boundary conditions.  For the spatial discretization,  we consider the mixed hybrid finite element method as proposed in~\cite{Radu11, Radu14}, in which the flux variable approximates the total flux (i.e. both diffusive and advective flux) and the Lagrange multiplier arising in the hybridization is used to discretize the advective term.  Such a mixed hybrid method is fully mass conservative,  as accurate as the standard mixed method~\cite{Dawson09} while it is more efficient in terms of computational cost and robust (with the use of an upwind operator) for problems with high Pecl\'et numbers.  
The semi-discrete interface problem with mixed hybrid finite elements is derived, and the fully discrete problem is obtained using backward Euler time-stepping.  The proposed method is fully implicit and global in time,  thus different time steps can be used in the subdomains; data will be exchanged from one time grid to another via a suitable $L^{2}$ projection in time.  Note that in this work we treat only conforming spatial discretization and focus on the use of local time stepping. The reader is referred to \cite{Mortar89, Mortar00, Mortar07, Mortar21}, where mortar mixed methods on nonmatching spatial grids are developed.  We shall investigate the accuracy and the convergence of the proposed method and compare the performance of optimized Ventcel and Robin conditions via numerical experiments with discontinuous coefficients and nonconforming time grids.  Analysis of the semi-discrete local Ventcel problem and the convergence of the iterative method is beyond the scope of this paper and will be considered separately in a forthcoming paper. 

The rest of the paper is organized as follows: after presenting the model problem of linear advection-diffusion equations,  we derive in Section~\ref{sec:Ventcellmixed} mixed formulations of Ventcel transmission conditions with nonoverlapping subdomains.  The global-in-time optimized Ventcel-Schwarz method at the continuous level is presented in Section~\ref{sec:contIP} and its space-discrete counterpart with mixed hybrid finite element discretization is discussed in Section~\ref{sec:semiIP}.  Nonconforming time discretization and the fully discrete interface problem are considered in Section~\ref{sec:time}. In Section~\ref{sec:NumRe},  numerical experiments are carried out to investigate the accuracy and convergence of the proposed method.  %Ventcel conditions and compare with Robin conditions on test cases with various Pecl\'et numbers and discontinuous coefficients. 

For a bounded domain $ \Omega $ of $ \mR^{d}$ $(d=2,3)$  with Lipschitz boundary $ \partial \Omega $ and some fixed time $T>0$,  we consider the following linear advection-diffusion problem
\begin{equation} \label{eq:model}
\begin{array}{rll} \omega \partial_{t} c + \nabla \cdot (\bu c - \bD \nabla c  ) & =f & \text{in} \; \Omega \times (0,T ),\\
c&=0 & \text{on} \; \partial \Omega\times (0,T), \\
c(\cdot, 0) & = c_{0} & \text{in} \; \Omega,
\end{array} 
\end{equation}
where $ c $ is the concentration of a contaminant dissolved in a fluid, $ f $ the source term, $ \omega $ the porosity, $ \bu $ the Darcy velocity ({\em assumed to be given and time independent}), $\bD$ a time-independent diffusion tensor.  We suppose that $\bD$ is diagonal and that each diagonal entry $D_{jj}$, $j = 1, . . . , d$, is
positive and bounded above and away from~$0$.  For simplicity,  we have imposed only Dirichlet boundary conditions; the formulations presented in the following can be generalized to other types of boundary conditions.  We rewrite \eqref{eq:model} in an equivalent mixed form by introducing the vector field $\bphi$, which consists of both diffusive and advective flux \cite{Radu11, Radu14}:
\begin{equation} \label{eq:mixed}
\begin{array}{rll} 
\omega \partial_{t} c + \nabla \cdot \bphi &=f  & \text{in} \; \Omega \times (0,T), \\
\bphi &= - \bD \nabla c + \bu c & \text{in} \; \Omega \times (0,T),
\end{array} 
\end{equation}
together with the boundary and initial conditions as in \eqref{eq:model}. 
%
%To define the corresponding weak formulation,  we denote by $(\cdot, \cdot)$ the inner product on $L^{2}(\Omega)$ or $(L^{2}(\Omega))^{2}$; for a measurable subset $\Theta \subset \Omega$, we write $(\cdot, \cdot)_{\Theta}$ and $\langle \cdot, \cdot \rangle_{\partial \Theta}$ to indicate the inner products considered on $\Theta$ and $\partial \Theta$ respectively.  
We denote by $(\cdot, \cdot)$ the inner product on $L^{2}(\Omega)$ or $(L^{2}(\Omega))^{2}$ and write the mixed variational formulation of~\eqref{eq:mixed} as follows:

For a.e. $t \in (0,T)$, find $\left (c(t), \bphi(t)\right ) \in L^{2}(\Omega)  \times H(\Div, \Omega)$ such that
\begin{equation} \label{eq:mixedweak}
\begin{array}{rll}
\left (\omega \partial_{t}c, \mu \right ) + \left (\nabla \cdot \bphi, \mu\right ) &= (f,\mu), & \forall \mu \in L^{2}(\Omega), \vspace{3pt}\\
\left (\bD^{-1} \bphi, \bv\right )-\left (\bD^{-1} \bu c, \bv\right ) - \left (c, \nabla \cdot \bv\right ) &=0, & \forall \bv \in H(\Div, \Omega).
\end{array} 
\end{equation}
For given $\bu \in (W^{1,\infty}(\Omega))^{2}$, $f \in C(0,T; L^{2}(\Omega))$ and $c \in H^{1}_{0}(\Omega)$,  there exists a unique solution to problem~\eqref{eq:mixedweak} as shown in \cite[Theorem 3.2]{Radu14}.  
\section{Domain decomposition with Ventcel transmission conditions in mixed form} \label{sec:Ventcellmixed}
\begin{figure}[!htbp]
\centering
\includegraphics[scale=0.9]{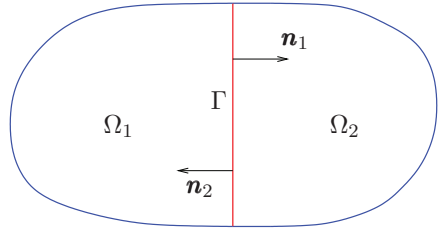}
\caption{Decomposition of $\Omega$ into two nonoverlapping subdomains.} 	
\label{Fig:NonoDD}  \vspace{-0.3cm}
\end{figure}

We consider a decomposition of $ \Omega $ into two nonoverlapping subdomains
$ \Omega_{1} $ and $ \Omega_{2} $ separated by an interface $\Gamma$ (see Figure~\ref{Fig:NonoDD}):
$$
\Omega_{1} \cap \Omega_{2} = \emptyset;\quad  \Gamma
= \partial \Omega_{1} \cap \partial \Omega_{2} \cap \Omega, \quad \Omega=  \Omega_{1} \cup \Omega_{2} \cup\Gamma. 
$$
The formulations given below can be generalized to the case of many subdomains in bands.

For $i=1,2$, let $ \pmb{n}_{i} $ denote the unit outward pointing normal vector field on $\partial\Omega_i$,
and for any scalar, vector or tensor valued function $v$ defined on $\Omega$, let $v_i$ be
the restriction  of $v$ to $ \Omega_{i} $.  
%Denote by $\bn_{\Gamma}=\bn_{1}\vert_{\Gamma}=-\bn_{2}\vert_{\Gamma}$ the normal unit vector to $\Gamma$ (see Figure~).
%We denote by $\tau_{j}$ $(j=1,\hdots, d-1)$ the mutually orthogonal unit tangential vectors to the interface $\Gamma$.  
In order to write the Ventcel transmission conditions, we use the
notation $ \nabla_{\tau}$ and $ \nabla_{\tau} \cdot \; $ for the tangential gradient
and divergence operators on~$\Gamma$ respectively.  For $i=1,2$, we denote by $\pmb{u}_{i,\Gamma}$ (and $\pmb{D}_{i,\Gamma}$)
%$\pmb{u}_{i,\Gamma}=\left (\pmb{u}_{i}-(\pmb{u}_{i}\cdot \bn_{i})\bn_{i}\right )\vert_{\Gamma}$ and $\pmb{D}_{i,\Gamma}=\left (\pmb{D}_{i}-(\pmb{D}_{i}\bn_{i})\otimes \bn_{i}\right )_{\Gamma}$
the tangential component of the trace of $\pmb{u}_{i}$ (respectively,  $\pmb{D}_{i} $) on $\Gamma$.  Problem~(\ref{eq:mixed}) can be reformulated as an equivalent multidomain problem consisting of the following space-time subdomain problems: 
\begin{equation} \label{eq:multi}
\begin{array}{rll} \omega_{i} \partial_{t} c_{i} + \nabla \cdot \bphi_{i}
        & =f & \text{in} \; \Omega_{i} \times (0,T),\\
\bphi_{i} &= - \bD_{i} \nabla c_{i} + \bu_{i} c_{i}  & \text{in} \; \Omega_{i} \times (0,T), \\
c_{i}    & = 0 & \text{on} \; \left (\partial \Omega_{i} \cap \partial \Omega \right) \times (0,T),\\
c_{i}(0) & = c_{0} & \text{in} \; \Omega_{i},
\end{array} \quad \text{for $ i=1,2 $}, \vspace{-0.1cm}
\end{equation}
together with the transmission conditions on the space-time interface: \vspace{-0.1cm}
\begin{equation} \label{eq:physicsTC}
\begin{array}{c}
c_{1} = c_{2} \\
\bphi_{1} \cdot \pmb{n}_{1} + \bphi_{2} \cdot \pmb{n}_{2} =0
\end{array} \quad \text{on} \; \Gamma \times \left (0,T\right ).  
\end{equation}

Under sufficient regularity, one may replace \eqref{eq:physicsTC} by the following Ventcel transmission conditions on $\Gamma \times \left (0,T\right )$:
\begin{equation}\label{eq:VentcelTCs}
\begin{array}{l}
-\bphi_{i} \cdot \bn_{i} + \alpha_{i,j} c_{i} + \beta_{i,j} \left (\omega_{j} \partial_{t} c_{i} +\nabla_{\tau} \cdot (\bu_{j,\Gamma} c_{i} - \bD_{j,\Gamma}\nabla_{\tau} c_{i})\right )\vspace{4pt} \\ \hspace{1cm}= -\bphi_{j} \cdot \bn_{i} + \alpha_{i,j} c_{j}+ \beta_{i,j} \left (\omega_{j} \partial_{t} c_{j} +\nabla_{\tau} \cdot (\bu_{j,\Gamma} c_{j} - \bD_{j,\Gamma}\nabla_{\tau}c_{j}) \right ),  
\end{array}
\end{equation}
for $i=1,2,$ and $j=(3-i)$, where $ \alpha_{i,j} $ and $ \beta_{i,j}$ are positive constants.  Ventcel or second-order transmission conditions in primal formulations were introduced for stationary problems in \cite{NR95,Japhet98} and then extended to time-dependent problems in~\cite{BGH09,HJO10,Japhet12}.  The parameters $ \alpha_{i,j} $ and $ \beta_{i,j}$ are chosen to optimize the convergence factor as studied in \cite{BGH09,HJO10}.  Note that the transmission conditions~\eqref{eq:VentcelTCs} reduce to Robin transmission conditions~\cite{H21} when $\beta_{i,j}=0$. 
%For the ease of presentation, we assume $\alpha_{i,j}=\alpha_{j,i}=\alpha$ and $\beta_{i,j}=\beta_{j,i}=\beta$; the performance of different choices of these coefficients shall be investigated in the numerical experiments.  
%To write the Ventcel conditions~\eqref{eq:VentcelTCs} fully in mixed form, we introduce the following new variables on the interface $ \Gamma $.  Let $ c_{i,\Gamma}, i=1,2, $ be the concentration trace $c_{i}$ on $ \Gamma $ and $\bphi_{i,\Gamma}:=\bu_{i,\Gamma} c_{i,\Gamma} - \bD_{i,\Gamma}\nabla_{\tau}c_{i,\Gamma}$ be the tangential component of the trace of $\bphi_{i}$ on $\Gamma$.  We especially define a vector field $\bphiG_{\Gamma, i}:= \bu_{j,\Gamma} c_{i,\Gamma}-\bD_{j, \Gamma} \nabla_{\tau} c_{i, \Gamma}, \, i=1,2, j=3-i. $ As the coefficients may be discontinuous across the interface,  $ \bphiG_{\Gamma, i} $ generally is not the tangential component of the trace of $ \bphi_{i} $ on the interface and it is used as an artificial tool for convergence purposes (it does not have a particular physical meaning).  We rewrite \eqref{eq:VentcelTCs} defined on $ \Gamma \times (0,T)$ as follows

To write the Ventcel conditions~\eqref{eq:VentcelTCs} fully in mixed form,  we introduce the following Lagrange multipliers on the interfaces $ \Gamma $: $ c_{i,\Gamma} $ representing the concentration trace $c_{i}$ on $ \Gamma $, and a vector field $\bphiG_{\Gamma, i}:= \bu_{j,\Gamma} c_{i,\Gamma}-\bD_{j, \Gamma} \nabla_{\tau} c_{i, \Gamma},$ for  $ i=1,2$ and $j=3-i. $ As the coefficients may be discontinuous across the interface,  $ \bphiG_{\Gamma, i} $ generally is not the tangential component of the trace of $ \bphi_{i} $ on the interface and it is used as an artificial tool for convergence purposes (it does not have a particular physical meaning).  Finally,  we denote by 
\begin{equation} \label{eq:tantrace}
\bpsi_{i,\Gamma}:=\left ( \bu_{i,\Gamma} - \bD_{i,\Gamma} \bD_{j,\Gamma}^{-1} \bu_{j,\Gamma} \right ) c_{i,\Gamma}+\bD_{i,\Gamma} \bD_{j,\Gamma}^{-1} \bphi_{\Gamma,i}, \; i=1,2; \, j=(3-i),
\end{equation}
the vector field representing the tangential component of the trace of $ \bphi_{i} $ on the interface.  We remark that $\bpsi_{i,\Gamma}$ is determined from $c_{i,\Gamma}$ and $\bphi_{\Gamma,i}$,  and it is used to exchange Ventcel data with the neighboring subdomain. 
%solving the subdomain problems does not involve $\bpsi_{i,\Gamma}$ (. 

The transmission conditions~\eqref{eq:VentcelTCs} on the space-time interface $\Gamma \times (0,T)$ can be rewritten as 
\begin{equation}\label{eq:VentcelTCsmixed}
\begin{array}{rl}
-\bphi_{i} \cdot \bn_{i} + \alpha_{i,j} c_{i,\Gamma} &\hspace{-0.2cm}+ \beta_{i,j} \left ( \omega_{j}\partial_{t} c_{i,\Gamma} +\nabla_{\tau} \cdot \bphiG_{\Gamma, i}\right ) \\
&= -\bphi_{j} \cdot \bn_{i} + \alpha_{i,j} c_{j,\Gamma} + \beta_{i,j} \left (\omega_{j} \partial_{t} c_{j,\Gamma} +\nabla_{\tau} \cdot \bpsi_{j,\Gamma} \right ),  \vspace{3pt} \\
\bphiG_{\Gamma, i}&=\bu_{j,\Gamma} c_{i,\Gamma} - \bD_{j,\Gamma}\nabla_{\tau} c_{i,\Gamma}, 
%\bpsi_{i,\Gamma}&=\left ( \bu_{i,\Gamma} - \bD_{i,\Gamma} \bD_{j,\Gamma}^{-1} \bu_{j,\Gamma} \right ) c_{i,\Gamma}+\bD_{i,\Gamma} \bD_{j,\Gamma}^{-1} \bphi_{\Gamma,i} ,
%\bu_{i,\Gamma} c_{i,\Gamma} - \bD_{i,\Gamma}\nabla_{\tau}c_{i,\Gamma},
\end{array}
\end{equation}
for $i=1,2,$, $j=(3-i)$. Next, we derive the formulation of the global-in-time optimized Ventcel-Schwarz method based on these transmission conditions.  We use the so-called Ventcel-to-Ventcel interface operators to rewrite the problems posed in the subdomains as a problem on the space-time interface.  

\section{Global-in-time optimized Ventcel-Schwarz method} \label{sec:contIP}
We introduce the interface unknowns
\begin{equation} \label{eq:zetai}
\zeta_{i}:=-\bphi_{i} \cdot \bn_{i} + \alpha_{i,j} c_{i} + \beta_{i,j} \left (\omega_{j} \partial_{t} c_{i} +\nabla_{\tau} \cdot \bphiG_{\Gamma, i}\right ) , \; \Gamma \times (0,T), \; \text{for} \; i=1,2.
\end{equation}
The transmission conditions~\eqref{eq:VentcelTCsmixed} become
\begin{equation}\label{eq:VentcelTCszeta}
\begin{array}{rl}
\zeta_{i} - \left [-\bphi_{j} \cdot \bn_{i} + \alpha_{i,j} c_{j}+ \beta_{i,j} \left (\omega_{j} \partial_{t} c_{j} +\nabla_{\tau} \cdot \bpsi_{j,\Gamma} \right )\right ]=0,   \; \text{on} \; \Gamma \times (0,T), 
\end{array}
\end{equation}
for $i=1,2,$ and $j=(3-i)$.  Equations~\eqref{eq:zetai} are used as Ventcel boundary conditions for the subdomain problems as presented in Subsection~\ref{subsec:locV}.  Then by enforcing the transmission conditions~\eqref{eq:VentcelTCszeta}, we obtain the space-time interface problem with the two unknowns $\zeta_{1}$ and $\zeta_{2}$ as discussed in Subsection~\ref{subsec:IP}.
\subsection{Local problem with Ventcel boundary conditions} \label{subsec:locV}
For a given function $\zeta \in L^{2}(0,T;\Theta_{\Gamma})$ with $\Theta_{\Gamma}:=L^{2}(\Gamma)$, consider the following advection-diffusion problem in subdomain $\Omega_{i}$ with Ventcel condition on the interface $\Gamma$:
\begin{equation} \label{eq:M2subVentcel}
\begin{array}{rll} \omega_{i}\partial_{t} c_{i} + \nabla \cdot \bphi_{i}
        & =f & \text{in} \; \Omega_{i} \times (0,T),\\
\bphi_{i} &= - \bD_{i} \nabla c_{i} + \bu_{i} c_{i}  & \text{in} \; \Omega_{i} \times (0,T), \\
c_{i}    & = 0 & \hspace{-1.4cm}\text{on} \; \left (\partial \Omega_{i} \cap \partial \Omega \right) \times (0,T),\\
-\bphi_{i} \cdot \bn_{i} + \alpha_{i,j} \; c_{i, \Gamma} &\hspace{-0.2cm}+ \beta_{i,j}  \left (\omega_{j} \partial_{t} c_{i,\Gamma}+ \nabla_{\tau} \cdot \bphiG_{\Gamma, i} \right ) = \zeta & \text{on} \; \Gamma \times (0,T), \\
\bphiG_{\Gamma , i}&= \bu_{j,\Gamma} c_{i,\Gamma}-\bD_{j, \Gamma} \nabla_{\tau} c_{i, \Gamma} & \text{on} \; \Gamma \times (0,T), \\
c_{i,\Gamma}&=0 & \text{on} \; \partial \Gamma \times (0,T), \\
c_{i}(0) & = c_{0} & \text{in} \; \Omega_{i},\\
c_{i,\Gamma}(0) & = c_{0}\vert_{\Gamma} & \text{on} \; \Gamma,
\end{array} 
\end{equation}
for $ i=1,2$ and $j=(3-i)$.  Problem~\eqref{eq:M2subVentcel} can be seen as a coupling of a $d$-dimensional PDE in the subdomain $\Omega_{i}$ and a $(d-1)$-dimensional PDE on the interface $\Gamma$; both PDEs are written in mixed form. 
To write the weak formulation for the local problem~\eqref{eq:M2subVentcel}, we introduce the following spaces:
\begin{align*}
 M_{i}&=\left \{ \overline{\mu}_{i} = (\mu_{i}, \mu_{i,\Gamma}) \in L^{2}(\Omega_{i}) \times L^{2}(\Gamma) \right \},\\
\Sigma_{i} & =\big \{  \overline{\pmb{v}}_{i} = (\pmb{v}_{i}, \pmb{v}_{\Gamma,i}) \in \pmb{L^{2}(\Omega_{i})} \times \pmb{L^{2}(\Gamma)}: \text{div}\; \pmb{v}_{i} \in L^{2} (\Omega_{i}) \; \, \text{and} \\
& \hspace{4.5cm} \; \; \beta_{i,j} \text{div}_{\tau} \; \pmb{v}_{\Gamma,i}- \pmb{v}_{i} \cdot \pmb{n}_{i}\vert_{\Gamma}  \in L^{2}(\Gamma) \big \}.
\end{align*}
For a measurable subset $W$ of $\Omega$,  we write $(\cdot, \cdot)_{W}$ to indicate the inner product on $W$. 
We define the following bilinear forms 
on $\Sigma_{i} \times \Sigma_{i} $, $\Sigma_{i} \times M_{i}$ and $M_{i} \times M_{i} $ respectively: \vspace{-0.2cm}
\begin{equation*}
\begin{array}{l}
	a_{i}(\overline{\bphi}_{i},\overline{\pmb{v}}_{i}) = \left (\pmb{D}_{i}^{-1} \bphi_{i}, \pmb{v}_{i}\right )_{\Omega_{i}} -\left (\bD_{i}^{-1}\bu_{i}c_{i}, \bv_{i}\right )_{\Omega_{i}}+ \left (\beta_{i,j} \bD_{j,\Gamma}^{-1}\bphiG_{\Gamma,i}, \pmb{v}_{\Gamma,i}\right )_{\Gamma} \\
	\hspace{6.5cm} -\left (\beta_{i,j} \bD_{j,\Gamma}^{-1}\bu_{j,\Gamma} c_{i,\Gamma},\pmb{v}_{\Gamma,i}\right )_{\Gamma},\vspace{3pt}\\
	 b_{i}(\overline{\bphi}_{i},\overline{\mu}_{i}) = \left (\nabla \cdot \bphi_{i},\mu_{i}\right )_{\Omega_{i}}+\left (\beta_{i,j} \nabla_{\tau} \cdot \bphiG_{\Gamma,i} - \bphi_{i} \cdot \bn_{i}, \mu_{i,\Gamma}\right )_{\Gamma}, \vspace{3pt}\\
	 \kappa_{i} (\overline{c}_{i}, \overline{\mu}_{i}) =\left (\omega_{i} c_{i}, \mu_{i}\right )_{\Omega_{i}}+ \left (\beta_{i,j} \omega_{j} c_{i, \Gamma}, \mu_{i,\Gamma}\right )_{\Gamma},  \quad \kappa_{i,\alpha}(\overline{c}_{i}, \overline{\mu}_{i}) = \left (\alpha_{i,j} c_{i, \Gamma}, \mu_{i,\Gamma}\right )_{\Gamma}, 
\end{array} 
\end{equation*}
and the linear form on $M_{i}$: 
\begin{equation*}
	\begin{array}{lccl}
	L(\overline{\mu}_{i})= \left (f, \mu_{i}\right )_{\Omega_{i}} + \left (\zeta,\mu_{i,\Gamma}\right )_{\Gamma}.
	\end{array} 
\end{equation*}
With the defined spaces and forms, the weak form of \eqref{eq:M2subVentcel} can be written as follows:\vspace{7pt}

Find $ \left (\overline{c}_{i}, \overline{\bphi}_{i}\right )  \in M_{i} \times \Sigma_{i} $ such that
\begin{eqnarray}
\begin{array}{rll} 
a_{i}(\overline{\bphi}_{i}, \overline{\pmb{v}}_{i}) - b_{i}(\overline{\pmb{v}}_{i}, \overline{c}_{i}) &=0 & \forall \overline{\pmb{v}_{i}} \in \Sigma_{i},\\
\kappa_{i}(\partial_{t}\overline{c}_{i}, \overline{\mu}_{i})+\kappa_{i,\alpha}(\overline{c}_{i} \overline{\mu}_{i}) +b_{i}(\overline{\bphi}_{i}, \overline{\mu}_{i}) & = L(\overline{\mu}_{i}) & \forall \overline{\mu}_{i} \in M_{i}.
\end{array}  \hspace{1cm} \label{eq:M2subVentcelweak}  
\end{eqnarray}
An existence and uniqueness result for evolution problems posed in mixed
form as~\eqref{eq:M2subVentcelweak} is analyzed in \cite{H16}.  In this paper,  we focus on finding the approximate solution using the mixed hybrid finite element method proposed in \cite{Radu11, Radu14} (see Section~\ref{sec:semiIP}).  
%$$ H^{1,1}_{\ast}(\Omega_{i}):=\left \{ \mu \in H^{1}(\Omega_{i}): \mu\vert_{\partial \Omega \cap \partial \Omega_{i}}=0 \; \text{and} \; \mu\vert_{\Gamma} \in H^{1}(\Gamma)\right \}. $$
\subsection{Space-time interface problem} \label{subsec:IP}
We aim to derive an interface problem associated with the subdomain problems~\eqref{eq:multi} with Ventcel transmission conditions~\eqref{eq:VentcelTCszeta}. Toward that end, we introduce the space
$$ H^{1,1}_{\ast}(\Omega_{i}):=\left \{ \mu \in H^{1}(\Omega_{i}): \mu\vert_{\partial \Omega \cap \partial \Omega_{i}}=0 \; \text{and} \; \mu\vert_{\Gamma} \in H^{1}(\Gamma)\right \},$$
and define the Ventcel-to-Ventcel operators $\iS^{\text{VtV}}_{i}$ as follows:
\begin{equation} 
\hspace{-0.3cm}\begin{array}{l} \iS^{\text{VtV}}_{i}: L^{2}(0,T;\Theta_{\Gamma}) \times L^{2}(0,T;L^{2}(\Omega_{i})) \times H^{1,1}_{\ast}(\Omega_{i})  \rightarrow L^{2}(0,T;\Theta_{\Gamma}) \vspace{0.1cm}\\
(\zeta, f,c_{0})  \longmapsto \mathcal{S}^{\text{VtV}}_{i}(\zeta, f,c_{0}) = -\bphi_{i} \cdot \pmb{n}_{j \mid \Gamma} + \alpha_{j,i} \; c_{i, \Gamma} + \beta_{j,i} \left ( \partial_{t} c_{i,\Gamma} +\nabla_{\tau} \cdot \bpsi_{i,\Gamma} \right ), 
\end{array}  \label{eq:S_VtV}
\end{equation}
where $\left (c_{i},\bphi_{i},c_{i,\Gamma}, \bphiG_{\Gamma,i}\right )$ is the solution of the local Ventcel problem \eqref{eq:M2subVentcel}, and $\bpsi_{i,\Gamma}$ is computed from $c_{i,\Gamma}$ and $\bphi_{\Gamma,i}$ as defined in \eqref{eq:tantrace}.
%Note that the solution $c_{i,\Gamma} \in H^{1}(\Gamma)$ as $$\nabla_{\tau} c_{i, \Gamma} = -\bD_{j, \Gamma}^{-1}\bphiG_{\Gamma , i}+\bD_{j, \Gamma}^{-1} \bu_{j,\Gamma} c_{i,\Gamma} \in L^{2}(\Gamma).$$
%
The interface problem, corresponding to the Ventcel transmission conditions~\eqref{eq:VentcelTCszeta}, is given by
\begin{equation}\label{eq:IP}
\begin{array}{ll}
\zeta_{1}=\iS^{\text{VtV}}_{2}(\zeta_{2},f,c_{0}), \\
\zeta_{2}=\iS^{\text{VtV}}_{1}(\zeta_{1},f,c_{0}),
\end{array} \quad \text{on} \; \Gamma \times (0,T).
\end{equation}

As the problems are linear, we can rewrite \eqref{eq:IP} equivalently as
\begin{equation}\label{eq:IPoperator}
\pmb{\mathcal{S}_{V}} \begin{pmatrix} \zeta_{1} \\ \zeta_{2} \end{pmatrix}  = \pmb{\mathcal{G}_{V}}(f,c_{0}),  \quad \text{on} \; \Gamma \times (0,T).
\end{equation}
where 
\begin{equation*} 
\begin{array}{cll} \pmb{\mathcal{S}}_{V}:\left (L^{2}(0,T;\Theta_{\Gamma})\right)^{2} & \longrightarrow &  \left (L^{2}(0,T;\Theta_{\Gamma})\right)^{2}\vspace{0.1cm}\\
\begin{pmatrix} \zeta_{1} \\ \zeta_{2} \end{pmatrix} & \longmapsto & \begin{pmatrix} \zeta_{1} -\mathcal{S}^{\text{VtV}}_{2}(\zeta_{2}, 0,0) \\
 \zeta_{2} -\mathcal{S}^{\text{VtV}}_{1}(\zeta_{1} , 0,0)
\end{pmatrix},
\end{array} 
\end{equation*}
and 
\begin{equation*} 
\begin{array}{cll} \pmb{\mathcal{G}}_{V}: L^{2}(0,T;L^{2}(\Omega)) \times H^{1,1}_{\ast}(\Omega)  & \longrightarrow &  \left (L^{2}(0,T;\Theta_{\Gamma})\right)^{2}\vspace{0.1cm}\\
(f,c_{0}) & \longmapsto & \begin{pmatrix} \mathcal{S}^{\text{VtV}}_{2}(0, f,c_{0}) \\
\mathcal{S}^{\text{VtV}}_{1}(0,f,c_{0}) 
\end{pmatrix},
\end{array} 
\end{equation*}
with $ H^{1,1}_{\ast}(\Omega):=\left \{ \mu \in H^{1}(\Omega): \mu\vert_{\partial \Omega}=0 \; \text{and} \; \mu\vert_{\Gamma} \in H^{1}(\Gamma)\right \}.$
The weak form of the space-time interface problem~\eqref{eq:IPoperator} is defined as
\begin{equation}\label{eq:IPoperatorweak}
\int_{0}^{T} \int_{\Gamma} \pmb{\mathcal{S}_{V}} \begin{pmatrix} \zeta_{1} \\ \zeta_{2} \end{pmatrix} \cdot \begin{pmatrix} \theta_{1} \\ \theta_{2} \end{pmatrix}  = \int_{0}^{T} \int_{\Gamma}  \pmb{\mathcal{G}_{V}}(f,c_{0})\cdot \begin{pmatrix} \theta_{1} \\ \theta_{2} \end{pmatrix} ,   \quad \forall (\theta_{1}, \theta_{2}) \in  \left (L^{2}(0,T;\Theta_{\Gamma})\right)^{2}. 
\end{equation}

%where 
%$$\bphiG_{\Gamma,i}=\bu_{i,\Gamma} c_{i} - \bD_{i,\Gamma}\nabla_{\tau}c_{i}.
%$$ 
%\textcolor{red}{For finite dimensional problems,  $\bphiG_{\Gamma,i}$ is computed via the weak form:
%$$\left (\bD_{i,\Gamma}^{-1}\bphiG_{\Gamma,i}, \bv_{\Gamma,i}\right )_{\Gamma}=\left (\bD_{i,\Gamma}^{-1}\bu_{i,\Gamma} c_{i}, \bv_{\Gamma,i}\right )_{\Gamma} + \left (\nabla_{\tau} \cdot \bv_{\Gamma,i}, c_{i}\right )_{\Gamma}.
%$$
%This gives a linear system in terms of the DOFs of $\bphiG_{\Gamma,i}$ (piecewise linear functions) - 1D mixed hybrid finite elements.
%}
%
The interface problem~\eqref{eq:IPoperatorweak} is solved iteratively by either Jacobi algorithm or a Krylov subspace iterative method (e.g., GMRES).  The former is equivalent to the OSWR algorithm with Ventcel conditions in mixed form, which reads as follows: starting with given initial guesses $g_{i,j} \in L^{2}(0,T;\Theta_{\Gamma})$ for the first iteration:
$$ g_{i,j}=-\bphi_{j}^{0} \cdot \bn_{i} + \alpha_{i,j} \; c_{j, \Gamma}^{0} + \beta_{i,j}  \left (\omega_{j} \partial_{t} c_{j,\Gamma}^{0}+ \nabla_{\tau} \cdot \bpsi_{j,\Gamma}^{0} \right ),
$$
then at each iteration $k=1,2, \hdots$, solve the subdomain problems for $i=1,2,$ and $j=3-i$:
\begin{equation} \label{eq:M2subVentcelOSWR}
\begin{array}{rll} \omega_{i}\partial_{t} c_{i}^{k} + \nabla \cdot \bphi_{i}^{k}
        & =f & \text{in} \; \Omega_{i} \times (0,T),\\
\bphi_{i}^{k} &= - \bD_{i} \nabla c_{i}^{k} + \bu_{i} c_{i}^{k}  & \text{in} \; \Omega_{i} \times (0,T), \\
c_{i}^{k}    & = 0 & \hspace{-1.6cm}\text{on} \; \left (\partial \Omega_{i} \cap \partial \Omega \right) \times (0,T),\\
-\bphi_{i}^{k} \cdot \bn_{i} + \alpha_{i,j} \; c_{i, \Gamma}^{k} &\hspace{-0.2cm}+ \beta_{i,j}  \left (\omega_{j} \partial_{t} c_{i,\Gamma}^{k}+ \nabla_{\tau} \cdot \bphiG_{\Gamma, i}^{k} \right ) & \\
& \hspace{-2cm}= -\bphi_{j}^{k-1} \cdot \bn_{i} + \alpha_{i,j} c_{j,\Gamma}^{k-1}+ \beta_{i,j} \left ( \omega_{j}\partial_{t} c_{j,\Gamma}^{k-1} +\nabla_{\tau} \cdot \bpsi_{j,\Gamma}^{k-1} \right ) & \text{on} \; \Gamma \times (0,T), \\
\bphiG_{\Gamma , i}^{k}&= \bu_{j,\Gamma} c_{i,\Gamma}^{k}-\bD_{j, \Gamma} \nabla_{\tau} c_{i, \Gamma}^{k} & \text{on} \; \Gamma \times (0,T), \\
c_{i,\Gamma}^{k}&=0 & \text{on} \; \partial \Gamma \times (0,T), \\
c_{i}^{k}(0) & = c_{0} & \text{in} \; \Omega_{i},\\
c_{i,\Gamma}^{k}(0) & = c_{0}\vert_{\Gamma} & \text{on} \; \Gamma,
\end{array} 
\end{equation}
where $\bpsi_{i,\Gamma}^{k-1}=\left ( \bu_{i,\Gamma} - \bD_{i,\Gamma} \bD_{j,\Gamma}^{-1} \bu_{j,\Gamma} \right ) c_{i,\Gamma}^{k-1}+\bD_{i,\Gamma} \bD_{j,\Gamma}^{-1} \bphi_{\Gamma,i}^{k-1}$.  The convergence of the Ventcel-OSWR algorithm as well as optimization of Ventcel parameters was studied in \cite{BGH09,Japhet12,BGGH16} with primal formulations.  We shall verify numerically the convergence of the iterates \eqref{eq:M2subVentcelOSWR} (after spatial and temporal discretizations) in Section~\ref{sec:NumRe}. 
%\bu_{i,\Gamma} c_{i,\Gamma} ^{k-1}- \bD_{i,\Gamma}\nabla_{\tau}c_{i,\Gamma}^{k-1}$.

% ------------------------------------------------
%
%    SECTION: Semi-discrete problem after spatial discretization by MHFE
%	
% ------------------------------------------------
%
\section{Space-discrete problems with mixed hybrid finite element discretization}
\label{sec:semiIP}
%
%
% DISCRETIZATION by MIXED HYBRID FEM
%
%
We consider the mixed hybrid finite element (MHFE) method proposed and analyzed in~\cite{Radu11,Radu14} for the spatial discretization of the local Ventcel problems~\eqref{eq:M2subVentcelweak}. The MHFE scheme is based on mixed finite elements together with the hybridization technique, in which the continuity constraint of the normal components of the fluxes over inter-element edges is relaxed via the use of Lagrange multipliers.  The Lagrange multipliers are also used to discretize the advection term, which is shown to be efficient for advection-dominant problems~\cite{Radu11,Radu14}.  Moreover, the MHFE formulation is very well-suited for using domain decomposition with Ventcel conditions since one has available both the normal trace of the flux and the trace of the concentration (i.e.  the Lagrange multiplier) on the interface.  

In the following, we consider the two-dimensional problem and assume $\Omega =\left (x_{a},x_{b}\right ) \times \left (y_{c},y_{d}\right )$ is a rectangular domain which is decomposed into two sub-rectangles $\Omega_{1}$ and $\Omega_{2}$ by a vertical interface $\Gamma=\{(x,y) \in \Omega: x = x_{\Gamma}\}$ for some $x_{\Gamma} \in \left (x_{a},x_{b}\right )$.  Note that the subdomain problems~\eqref{eq:M2subVentcelweak} with Ventcel conditions require MHFE discretization in both one and two dimensions.
%We assume $\Omega$ is a rectangle in $\mR^{2}$ %divided into two subdomains with a vertical interface.
Let $ \iK_{h,i}, \; i=1,2, $ be a finite element partition of each $ \Omega_{i} $ into rectangles such that their
union $ \iK_{h}= \cup_{i=1}^2 \iK_{h,i}$ forms a finite element partition of $\Omega $.  We assume the spatial discretization is conforming and denote by $ \iE_{h}^{\Gamma} $ the set of edges of elements of $\iK_{h,1}$ or $\iK_{h,2}$ that lie on $ \Gamma$. Let $\iE_{h,i}$ be the set of all edges of elements of $ \iK_{h,i}$:
$$\iE_{h,i}=\iE_{h,i}^{I} \cup \iE_{h,i}^{D} \cup \iE_{h}^{\Gamma},$$
where $\iE_{h,i}^{I}$ the set of all interior edges and $\iE_{h,i}^{D}$ the set of edges on the external boundary $\partial \Omega_{i} \cap \partial \Omega$, for $i=1,2$.  Denote by $\iP_{h}^{\Gamma}$ the set of endpoints $P$ of interface edges $E \in \iE_{h}^{\Gamma}$.  
For $K \in \iK_{h} $, let $\bn_K$ be the unit, normal, outward-pointing vector field on the boundary $\partial K$; for each edge $E \subset \partial K$, we denote by $\bn_{E}$ the unit normal vector of $E$, outward to $K$. Let $h_{K}=\text{diam}(K)$ and $h=\max_{K \in \iK_{h}} h_{K}$. 
The discrete spaces for the two-dimensional scalar and vector variables are defined based on the lowest-order Raviart-Thomas mixed finite elements:
\begin{align*}
M_{h,i}&:=\left \{ \mu \in L^{2}(\Omega_{i}): \mu_{\mid K}=\text{constant}, \; \forall K \in \iK_{h,i} \right \}, \\
\Sigma_{h,i} & := \left \{ \bv \in (L^{2}(\Omega_{i}))^{2}: \bv\vert_{K} \in \Sigma_{K} , \; \forall K \in \iK_{h,i} \right \},
\end{align*}
where $\Sigma_{K}:= \left \{ \bv:  K \rightarrow \mR^{2}, \; \bv = \left (a_{K} + b_{K} x, a^{\prime}_{K} +b^{\prime}_{K} y \right ), \, (a_{K}, b_{K}, a^{\prime}_{K}, b^{\prime}_{K}) \in \mR^{4} \right \}$, for $K \in \iK_{h,i},$ is the local Raviart-Thomas space. 
Note that for the sake of presentation simplicity,  we have omitted the use of the subscript $h$ for functions of the finite dimensional spaces.  The discrete space for the Lagrange multiplier representing the trace of the concentration on the edges is given by
$$
\Lambda_{h,i} := \left \{ \eta \in L^{2}(\iE_{h,i}): \eta\vert_{E}=\, \text{constant on E}, \; \forall E \in \iE_{h,i} \; \text{and } \, \eta\vert_{E}=0, \; \forall E \in \iE_{h,i}^{D} \right \}.
$$
In addition,  to take into account the interface as part of the subdomain boundary we define the space
$$ \Lambda_{h,i}^{\Gamma,0} := \left \{ \eta \in  \Lambda_{h,i}:  \eta\vert_{E}=0, \; \forall E \in \iE_{h}^{\Gamma}\right \}. 
$$
For the interface, the one-dimensional Raviart-Thomas spaces are given by
\begin{align*}
\Theta_{h, \Gamma} &= \left \{ \theta \in L^{2}(\Gamma): \theta\vert_{E}=\, \text{constant on E}, \; \forall E \in \iE_{h}^{\Gamma}  \right \}, \\
\Sigma_{h,\Gamma} & := \left \{ \bv_{\Gamma} \in L^{2}(\Gamma): \bv_{\Gamma}\vert_{E} \in \Sigma_{\Gamma,E} , \; \forall E \in \iE_{h}^{\Gamma} \right \},
\end{align*}
where $\Sigma_{\Gamma,E}:= \left \{ \bv_{\Gamma}:  E \rightarrow \mR, \; \bv_{\Gamma} = a_{E} + b_{E}y, \, (a_{E}, b_{E}) \in \mR^{2} \right \},$ for $E \in \iE_{h}^{\Gamma}, $ is the local one-dimensional Raviart-Thomas space.  Note that we still use bold fonts for the one-dimensional flux variable and test functions to be consistent with the formulations given in the previous section.
The space of the Lagrange multiplier for the one dimensional equation is defined as:
$$ \Xi_{h,\Gamma}=\left \{ \varsigma: \iP_{h}^{\Gamma} \rightarrow \mR, \; \, \varsigma(P)=0 \; \text{if}\; P \in \partial \Gamma \right \}.
$$

For $ c_{i}(t) \in M_{h,i} $ and $\lambda_{i}(t) \in \Lambda_{h,i}$,  we have the representation 
$$c_{i} (t,x,y)= \sum_{K \in \iK_{h,i}} c_{i,K}(t) \chi_{K}(x,y), \quad \lambda_{i} (t,y)= \sum_{E \in \iE_{h,i}} \lambda_{i,E}(t) \chi_{E}(y),$$ 
where $\chi_{K}$ and $\chi_{E}$ are the characteristic functions of element $K \in \iK_{h,i}$ and edge $E \in \iE_{h,i}$, respectively.  The values $c_{i,K}$ and  $\lambda_{i,E}$ represent the average of the concentration on $K$ and on $ E $,  respectively. 

For $\bphi_{i}(t) \in \Sigma_{h,i}$, the function is defined locally as
$$ \bphi_{i} (t,x,y) \vert_{K}= \sum_{E \subset \partial K} \varphi_{i,KE} (t)\bw_{KE}(x,y), 
$$
where $\varphi_{i,KE}$ is the normal flux leaving $K$ through the edge $E$ and $\{\bw_{KE}\}_{E \subset \partial K}$ are the basis functions of the local Raviart-Thomas space $\Sigma_{K}$ satisfying
$$ \int_{E^{\prime}} \bw_{KE} \cdot \bn_{K} = \delta_{E,E^{\prime}}, \; \forall E^{\prime} \subset \partial K. 
$$

Similarly, for $\bphiG_{\Gamma,i}(t) \in \Sigma_{h,\Gamma} $,  we have the expression 
$$\bphiG_{\Gamma,i} (t,y)= \sum_{P \in \partial E} \varphi_{\Gamma,i,EP} (t)\; \bw_{\Gamma,EP}(y),
$$ 
where $\{\bw_{\Gamma,EP}\}_{P \in \partial E}$ are the basis functions of $\Sigma_{\Gamma,E}$.
We still denote by $\bu_{i}$ the projection of the Darcy velocity on $\Sigma_{h,i}$
% which is defined as 
%  u_{KE} is the flow rate of $\bu$ through $E$
and by $\bu_{i,\Gamma}$ the projection of its tangential component on $\Sigma_{h,\Gamma}$:
$$\bu_{i}:=\sum_{K \in \iK_{h,i}} \sum_{E \subset \partial K} u_{i,KE}\bw_{KE}(\bx), \quad \bu_{i,\Gamma}:=\sum_{E \in \iE^{\Gamma}_{h}} \sum_{P \in \partial E}u_{i,\Gamma,EP}  \; w_{\Gamma,EP}(y).$$
\subsection{Semi-discrete local Ventcel problems} \label{subsec:semilocV}
The mixed hybrid formulation for the local Ventcel problem~\eqref{eq:M2subVentcelweak} is given by: 

Find $\left (c_{i}(t), \bphi_{i}(t), \lambda_{i}(t),\bphiG_{\Gamma,i}(t), \xi_{\Gamma,i}(t)\right ) \in M_{h,i} \times \Sigma_{h,i} \times \Lambda_{h,i} \times \Sigma_{h,\Gamma} \times \Xi_{h,\Gamma}$, for a.e. $t \in (0,T)$,  such that
\begin{align}
&\left (\omega_{i} \partial_{t}c_{i}, \mu\right )_{\Omega_{i}} + \left (\nabla \cdot \bphi_{i}, \mu\right )_{\Omega_{i}} = (f,\mu)_{\Omega_{i}}, \quad \forall \mu \in M_{h,i}, \vspace{3pt} \label{eq:semiVmass}\\
&\left (\bD_{i}^{-1} \bphi_{i}, \bv\right )_{\Omega_{i}}-\sum_{K \in \iT_{h,i}} \sum_{E \subset \partial K} u_{i,KE} \; \lambda_{i,E} \left (\bD_{i}^{-1} \bw_{KE}, \bv\right )_{K} - \left (c_{i}, \nabla \cdot \bv\right )_{\Omega_{i}} \nonumber\\
& \hspace{4cm} + \sum_{K \in \iK_{h}} \left (\lambda_{i}, \bv \cdot \bn_{K} \right)_{\partial K} =0, \quad \forall \bv \in \Sigma_{h,i}, \vspace{3pt} \label{eq:semiVflux}\\
&\sum_{K \in \iK_{h,i}} \left ( \eta, \bphi_{i} \cdot \bn_{K}\right)_{\partial K} = 0, \quad \forall \eta \in \Lambda_{h,i}^{\Gamma,0},\vspace{3pt} \label{eq:semiVlag}\\
&\left (-\bphi_{i} \cdot \bn_{i}, \theta \right )_{\Gamma}+ \left (\alpha_{i,j} \; \lambda_{i}, \theta \right )_{\Gamma} + \left (\beta_{i,j}  \omega_{j} \partial_{t} \lambda_{i}, \theta \right )_{\Gamma} +\left (\beta_{i,j} \nabla_{\tau} \cdot \bphiG_{\Gamma, i},\theta  \right )_{\Gamma}  \nonumber \vspace{3pt}\\
&\hspace{6.5cm} = \left (\zeta,\theta \right )_{\Gamma}, \; \forall \theta  \in \Theta_{h,\Gamma}, \label{eq:semiV1dmass}\vspace{3pt}\\
& \left (\bD_{j, \Gamma} ^{-1}\bphiG_{\Gamma , i}, \bv_{\Gamma}\right )_{\Gamma} -\sum_{E \in \iE_{h}^{\Gamma}} \sum_{P \in \partial E} u_{i,\Gamma,EP} \; \xi_{i,P} \left (\bD_{j,\Gamma}^{-1} \bw_{\Gamma,EP}, \bv_{\Gamma}\right )_{E}  \nonumber \\
&\hspace{2cm}- \left (\lambda_{i}, \nabla_{\tau} \cdot \bv_{\Gamma}\right )_{\Gamma} + \sum_{E \in \iE_{h}^{\Gamma}} \left (\xi_{i}\, \bv_{\Gamma}\right )\vert_{\partial E} =0, \quad \forall \bv_{\Gamma} \in \Sigma_{h,\Gamma}, \vspace{3pt} \label{eq:semiV1dflux}\\
&\sum_{E \in \iE_{h}^{\Gamma}} \left (\varsigma\, \bphiG_{\Gamma,i}\right ) \vert_{\partial E}  = 0, \quad \forall \varsigma \in \Xi_{h}^{\Gamma}. \label{eq:semiV1dlag}
\end{align}

Note that for the space-discrete advection term in \eqref{eq:semiVflux}, we have used the Lagrange multiplier instead of using the piecewise constant concentration. Such a scheme is shown to give good numerical performance for the case where advection is moderately dominant \cite{Radu11, H21}.  For strongly advection-dominated problems, using upwind values are recommended; interested readers are referred to \cite{Radu11, Radu14} for further details. 
In addition, equation~\eqref{eq:semiVlag} enforces the continuity of the normal components of the fluxes over inter-element edges so that the vector variable $ \bphi_{i} \in \Sigma_{h,i}$ belongs to $H(\text{div}, \Omega_{i})$.  Similar treatment is done for the one-dimensional equations on the interface (cf.  \eqref{eq:semiV1dmass}-\eqref{eq:semiV1dlag}).
By taking the test functions to be basis functions in \eqref{eq:semiVmass}-\eqref{eq:semiV1dlag}, we obtain a system of linear equations as shown in Appendix~\ref{app:algsystem}.  Next we formulate the semi-discrete interface problem with MHFE discretization.

\subsection{Semi-discrete continuous-in-time interface problem}
%The semi-discrete counterpart of the transmission conditions~\eqref{eq:VentcelTCszeta} is given by
%\begin{equation} \label{eq:semiTCs}
%\begin{array}{l}
%(\zeta_{i},\theta)_{\Gamma}-  \big [\left (-\bphi_{j} \cdot \bn_{i}, \theta \right )_{\Gamma}+ \left (\alpha_{i,j} \; \lambda_{j}, \theta \right )_{\Gamma} + \left (\beta_{i,j}  \omega_{j} \partial_{t} \lambda_{j}, \theta \right )_{\Gamma} \vspace{3pt}\\
%\hspace{4cm}+\left (\beta_{i,j} \nabla_{\tau} \cdot \varphi_{j,\Gamma},\theta  \right )_{\Gamma} \big]=0,  \quad  \forall \theta \in \Theta_{h,\Gamma},
% \end{array}
%\end{equation}
%\begin{equation} \label{eq:semiTCs}
%(\zeta_{i},\theta)_{\Gamma}-  \left (-\bphi_{j} \cdot \bn_{i} +\alpha_{i,j} \; \lambda_{j} +\beta_{i,j}  \left ( \omega_{j} \partial_{t} \lambda_{j}+ \nabla_{\tau} \cdot \varphi_{j,\Gamma}\right ),\theta  \right )_{\Gamma} =0,  \quad  \forall \theta \in \Theta_{h,\Gamma},
%\end{equation}
%for $i=1,2$, $j=(3-i)$ and for a.e.  $t \in (0,T)$, where $\left (c_{i}(t), \bphi_{i}(t), \lambda_{i}(t),\bphiG_{\Gamma,i}(t), \xi_{i}(t)\right ) \in M_{h,i} \times \Sigma_{h,i} \times \Lambda_{h,i} \times \Sigma_{h,\Gamma} \times \Xi_{h,\Gamma}$ is the solution to \eqref{eq:semiVmass}-\eqref{eq:semiV1dlag} with $\zeta = \zeta_{i}$, for $i=1,2$.
The semi-discrete Ventcel-to-Ventcel operators $\iS_{h,i}^{\text{VtV}}$ are given by
\begin{equation}\label{eq:semiSvtv}
 \iS_{h,i}^{\text{VtV}}(\zeta,f,c_{0})=-\bphi_{i} \cdot \bn_{j}+ \alpha_{j,i} \; \lambda_{i}+ \beta_{j,i} \left ( \omega_{i} \partial_{t} \lambda_{i}+ \nabla_{\tau} \cdot \bpsi_{i,\Gamma}\right ) \in L^{2}(0,T;\Theta_{h,\Gamma}),
\end{equation}
where $\left (c_{i}(t), \bphi_{i}(t), \lambda_{i}(t),\bphiG_{\Gamma,i}(t), \xi_{\Gamma, i}(t)\right ) \in M_{h,i} \times \Sigma_{h,i} \times \Lambda_{h,i} \times \Sigma_{h,\Gamma} \times \Xi_{h,\Gamma}$ is the solution to \eqref{eq:semiVmass}-\eqref{eq:semiV1dlag}, for $i=1,2$, and $ \bpsi_{i,\Gamma} \in \Sigma_{h,\Gamma}$ is determined from $\lambda_{i}, \bphiG_{\Gamma,i}$ and $\xi_{\Gamma, i}$ by
\begin{align} 
 &\left (\bpsi_{i,\Gamma}, \bv_{\Gamma}\right )_{\Omega_{i}}= \sum_{E \in \iE_{h}^{\Gamma}} \sum_{P \in \partial E} u_{i,\Gamma,EP} \; \xi_{i,P} \left (\bw_{\Gamma,EP}, \bv_{\Gamma}\right )_{E} \nonumber\\
& - \sum_{E \in \iE_{h}^{\Gamma}} \sum_{P \in \partial E} u_{j,\Gamma,EP} \; \xi_{i,P} \left (\bD_{i,\Gamma}\bD_{j,\Gamma}^{-1} \bw_{\Gamma,EP}, \bv_{\Gamma}\right )_{E} +\left (\bD_{i,\Gamma} \bD_{j,\Gamma}^{-1} \bphi_{\Gamma,i}, \bv_{\Gamma}\right )_{\Omega_{i}},\label{eq:tantracedis}
\end{align}
for all $\bv_{\Gamma} \in \Sigma_{h,\Gamma}$.  It should be noted that $\bpsi_{i,\Gamma}$ is not computed explicitly,  i.e. we do not solve \eqref{eq:tantracedis}; instead we will use \eqref{eq:tantracedis} to calculate the Ventcel data~\eqref{eq:semiSvtv} as detailed in the Appendix~\ref{app:comphi}.

The space-discrete counterpart of the interface problem~\eqref{eq:IPoperatorweak} is as follows:  for a.e.  $t \in (0,T)$,
\begin{equation} \label{eq:semiIP}
 \begin{array}{rl} 
 (\zeta_{1}, \theta_{1})_{\Gamma} - (\iS_{h,2}^{\text{VtV}}(\zeta_{2},0,0),\theta_{1})_{\Gamma}&=(\iS_{h,2}^{\text{VtV}}(0,f,c_0),\theta_{1})_{\Gamma}, \vspace{4pt}\\ 
 (\zeta_{2}, \theta_{2})_{\Gamma} - (\iS_{h,1}^{\text{VtV}}(\zeta_{1},0,0),\theta_{2})_{\Gamma} &=(\iS_{h,1}^{\text{VtV}}(0,f,c_0),\theta_{2})_{\Gamma},
 \end{array} \quad \forall (\theta_{1}, \theta_{2}) \in \left (\Theta_{h,\Gamma}\right )^{2}.
\end{equation}
%or equivalently 
%\begin{equation}\label{eq:semiIP}
%\int_{0}^{T} \int_{\Gamma} \pmb{\mathcal{S}_{h,V}} \begin{pmatrix} \zeta_{1} \\ \zeta_{2} \end{pmatrix} \cdot \begin{pmatrix} \theta_{1} \\ \theta_{2} \end{pmatrix}  = \int_{0}^{T} \int_{\Gamma}  \pmb{\chi_{h,V}}(f,c_{0})\cdot \begin{pmatrix} \theta_{1} \\ \theta_{2} \end{pmatrix} ,   \quad \forall (\theta_{1}, \theta_{2}) \in \Theta_{h,\Gamma} \times \Theta_{h,\Gamma}. 
%\end{equation}
Again,  we solve this interface problem iteratively by Jacobi or GMRES; at each iteration,  the semi-discrete local Ventcel problem \eqref{eq:semiVmass}-\eqref{eq:semiV1dlag} is solved over the whole time interval in each subdomain.  In the next section, we consider the fully discrete interface problem with nonconforming time grids.
%
% ------------------------------------------------
%
%    SECTION: Nonconforming time grids
%	
% ------------------------------------------------
%
\section{Nonconforming time discretizations}
\label{sec:time}
\begin{figure}[h]
\vspace{1.2cm}
\centering
\begin{minipage}[c]{0.5 \linewidth}
\setlength{\unitlength}{1pt} 
\begin{picture}(140,70)(0,0)
\thicklines
\put(0,3){\includegraphics[scale=0.55]{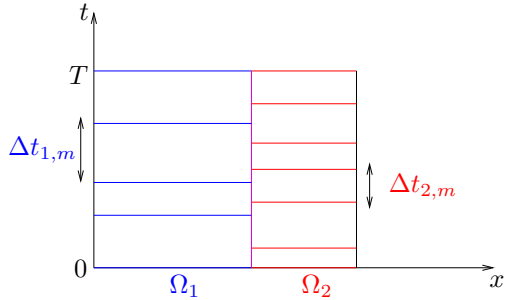} \\}
\put(-1,1){$ 0 $}
\put(-3,74){$ T $}
\put(35,-5){\textcolor{blue}{$ \Omega_{1} $}}
\put(85,-5){\textcolor{red}{$ \Omega_{2}$}}
\put(-26,47){\textcolor{blue}{$\Delta t_{1,m} $}}
\put(118,33){\textcolor{red}{$ \Delta t_{2,m} $}}
\put(156,-3){$ x $}
\put(1,98){$ t $}
\end{picture}
\end{minipage}
\caption{Nonconforming time grids in the subdomains.}
\label{Fig:Time} \vspace{-0.4cm}
\end{figure}

As the interface problem~\eqref{eq:semiIP} is global in time, independent time discretizations can be used in the subdomains.  Let $ \iT_{1} $ and $ \iT_{2} $
be two different partitions of the time interval $ (0,T) $ into sub-intervals (see Figure~\ref{Fig:Time}).
We denote by $ J_{i,m} $ the time interval $ (t_{i,m}, t_{i,m-1}] $ and by
$ \Delta t_{i,m} := (t_{i,m} - t_{i,m-1}) $ for $ m=1, \hdots, M_{i} $ and $ i=1,2 $. 
We use the backward Euler method to advance in time implicitly; the same idea can be generalized to higher order methods. 
\noindent We denote by $ P_{0}(\mathcal{T}_{i}, \Theta_{h,\Gamma}) $ the space of piecewise constant functions in time on grid
$ \mathcal{T}_{i} $ with values in $ \Theta_{h,\Gamma} $.
%\begin{equation} \label{P0space}
%P_{0}(\mathcal{T}_{i}, W) = \left \{ \psi: (0,T) \rightarrow W,
%\psi \; \; \text{ is constant on} \;  J_{m}^{i}, \ \forall m=1, \dots, M_{i} \right \}. 
%\end{equation}
To exchange data on the space-time interface between different time grids, we use an average-valued projection $ \Pi_{ji} $ from  $ P_{0} (\mathcal{T}_{i}, \Theta_{h,\Gamma}) $ onto $ P_{0}(\mathcal{T}_{j},\Theta_{h,\Gamma}) $ (see \cite{OSWRwave03,Japhet12}):
%~: for $ \psi \in P_{0} (\mathcal{T}_{i}, W)$,$ \Pi_{ji} \psi \hspace{-2pt} \mid_{J^{j}_{m}} $ is the average value of $ \psi $ on $ J^{j}_{m} $,for $ m=1, \dots, M_{j} $: 
%
\begin{equation} \label{ProjectionTime}
\Pi_{ji} \left ( \psi \right )\mid_{J_{j,m}}
  = \frac{1}{\Delta t_{j,m}} \sum_{l=1}^{M_{i}} \int_{J_{j,m} \cap J_{i,l}} \psi, \quad \text{for} \; \psi \in P_{0} (\mathcal{T}_{i}, \Theta_{h,\Gamma}).
\end{equation}
This projection is implemented numerically using the algorithm given in \cite{Projection1d}.
%With these tools, we are now able to weakly enforce the transmission conditions over the time intervals. We still denote by $ (c_{i}, \br_{i}) $, for $ i=1,2 $, the solution of the problem semi-discrete in time corresponding to problem~\eqref{Schursub-mixed} or~\eqref{Schwarzsub-mixed}. 

Let $\zeta_{i} =  (\zeta_{i}^{m})_{m=1}^{M_{i}} \in P_{0}(\mathcal{T}_{i}, \Theta_{h,\Gamma})$, for $i=1,2$. The fully discrete counterpart of the interface problem~\eqref{eq:semiIP} is obtained by weakly enforcing the transmission conditions over the time intervals:
%\begin{align}
%\int_{J_{1,m}} \left [(\zeta_{1}, \theta_{1})_{\Gamma} - (\Pi_{12}\left (\iS_{h,2}^{\text{VtV}}(\zeta_{2},0,0)\right ),\theta_{1})_{\Gamma}\right ] dt&=\int_{J_{1,m}}(\Pi_{12}\left (\iS_{h,2}^{\text{VtV}}(0,f,c_0)\right ),\theta_{1})_{\Gamma} dt.
%\end{align}
\begin{equation}\label{eq:fulldisIP}
\begin{array}{ll}
&\hspace{-1cm}\int_{J_{i,m}} \left [(\zeta_{i}, \theta_{i})_{\Gamma} - \left (\Pi_{ij}\left (\iS_{h,j}^{\text{VtV}}(\zeta_{j},0,0)\right ),\theta_{i}\right )_{\Gamma}\right ] dt \vspace{3pt}\\
&=\int_{J_{i,m}}\left (\Pi_{ij}\left (\iS_{h,j}^{\text{VtV}}(0,f,c_0)\right ),\theta_{i}\right )_{\Gamma} dt, \quad \forall \theta_{i} \in \Theta_{h,\Gamma}, \forall m=1, \hdots, M_{i}, 
\end{array}
\end{equation}
where $\iS_{h,i}^{\text{VtV}}(\zeta_{i},f,c_{0}) \in P_{0}(\mathcal{T}_{i}, \Theta_{h,\Gamma})$ is computed from the solution to the fully discrete subdomain problem obtained by applying the backward Euler method to \eqref{eq:semiV1dmass}-\eqref{eq:semiV1dlag} on time grid $\iT_{i}$:
$$\iS_{h,i}^{\text{VtV}}(\zeta_{i},f,c_{0})=\left (-\bphi_{i}^{m} \cdot \bn_{j}+ \alpha_{j,i} \; \lambda_{i}^{m} + \beta_{j,i} \left ( \omega_{i} \left (\frac{\lambda_{i}^{m}-\lambda_{i}^{m-1}}{\Delta t^{i,m}}\right )+ \nabla_{\tau} \cdot \bpsi_{i,\Gamma}^{m}\right )\right )_{m=1}^{M_{i}},
$$ 
for $i=1,2$ and $j=(3-i)$.

%
% ------------------------------------------------
%
%    SECTION: Numerical results
%	
% ------------------------------------------------
%
\section{Numerical experiments}
\label{sec:NumRe}
We study the numerical performance of the global-in-time optimized Ventcel-Schwarz method presented in previous sections.  In our numerical experiments, we consider the diffusion tensor $ \bD_{i} = d_{i} \pmb{I} $ isotropic and constant on each subdomain, where $ \pmb{I} $ is the 2D identity matrix.  The domain of calculation $\Omega$ is the unit square which is decomposed into two equal subdomains $ \Omega_{1} = (0,0.5) \times (0,1) $ and $ \Omega_{2} = (0.5,1) \times (0,1) $.  For the spatial discretization, we consider a conforming rectangular mesh with size $ h_{1}=h_{2}= h$.  For the time discretization,  we use different time steps $ \Delta t_{1} \neq \Delta t_{2}$.  Two test cases are presented: Test case~1 with the same constant coefficients in the subdomains to verify the accuracy of the method,  and Test case~2 with various discontinuous coefficients corresponding to different Pecl\'et numbers to check the robustness of the method when advection is dominant. 

For Ventcel parameters,  we consider two choices: the optimized {\em one-sided} Ventcel parameters, i.e. $\alpha_{i,j}=p^{\ast}$ and $\beta_{i,j}=q^{\ast}$,  and the optimized {\em weighted} Ventcel parameters given by:
\begin{equation}
\alpha_{i,j}=p^{\ast}-\frac{\bu_{j} \cdot \bn_{j}}{2},  \quad \beta_{i,j}=d_{j} q^{\ast}, \quad \text{for} \; i=1,2, \; j=(3-i), 
\end{equation}
where $p^{\ast}$ and $q^{\ast}$ are positive numbers and are obtained by numerical optimization of the theoretical convergence factor \cite{GHK07,HJO10}.  We also include the results by the optimized two-sided Robin transmission conditions in \cite{H21}, i.e.  with $\alpha_{i,j} \neq \alpha_{j,i}$ and $\beta_{i,j}=0$, for comparison purposes.
%
%
%	 TEST: analytical solution
%
%
	\subsection{Test case 1: with a known analytical solution}
	
%
% ----------- 1. Description ----------------
%	
We first verify the accuracy in space and in time of the proposed algorithm by considering a test case with the exact solution is given by $$ u(x,y,t)=e^{-4t} \, \sin (\pi x) \, \sin (\pi y).$$%Dirichlet conditions are imposed on the boundary of $\Omega$.
Constant parameters are imposed on the whole domain: $ \omega_{i}=1 $, $ \bu_{i} =(1, \; 1)^{T} $,  and $d_{i}=1$, for $i=1,2$.  The interface problem \eqref{eq:fulldisIP} is solved iteratively by either Jacobi iteration or GMRES with a zero initial guess on the interface; the iteration stops when the relative residual is smaller than $10^{-6}$. 
%
% Convergence in space
%
In Table~\ref{tab:test1accuSpace}, we show the relative $L^{2}(\Omega)-$norm errors of $c$ and $\bphi$ at $T=0.1$ with fixed time step sizes $\Delta t_{1}=T/80$ and $\Delta t_{2}=T/60$ and a decreasing mesh size $h$.  These errors are obtained using the optimized Ventcel transmission conditions.  We observe that the order of accuracy in space is preserved with nonconforming time grids.  The numbers of Jacobi and GMRES iterations are also reported in Table~\ref{tab:test1accuSpace} with different choices of the optimized parameters.  The convergence of OSWR methods with optimized Ventcel parameters is almost independent of the mesh size. This is consistent with the theoretical result in \cite{BGH09}, where it is shown that the convergence factor of the Ventcel-OSWR method behaves like $1-O(h^{1/4})$.  Moreover,  the Ventcel-OSWR method converges much faster, by nearly a factor of 2,  than the Robin-OSWR method.  As the coefficients are continuous,  there is not much difference in terms of convergence speed between the optimized one-sided and weighted Ventcel parameters.  GMRES improves slightly the convergence compared to Jacobi iteration. 
%
% PURE DIRICHLET BCs
%
%\begin{table}[h] 
%\centering
%\setlength{\extrarowheight}{2pt}
%	\begin{tabular}{| l | c | c | c | c | } \hline 
%		\multirow{2}{*}{$h$} & \multicolumn{2}{c|}{$L^{2}$ errors} & \multicolumn{2}{c|}{Number of iterations}\\\cline{2-3} \cline{4-5}
%		& $c$ & $\bphi$ & Jacobi & GMRES \\ \hline
%		%
%		$1/20$ &  0.0641	\phantom{[1.00]}&  0.0453	\phantom{[1.00]} & 11  & 10  \\
%		$1/40$ & 0.0321	[1.00] & 0.0227 [1.00]& 11 & 11\\
%		$1/80$ &  0.0160	[1.00] & 0.0114 [0.99]& 12 &11\\
%		$1/160$ & 0.0080	[1.00] & 0.0057 [1.00]& 13 &12 \\ \hline
%	\end{tabular}  
%	\caption{[Test case 1] Accuracy in space, the convergence rates are shown in square brackets.} 
%\label{tab:test1accuSpace}
%\end{table}

\begin{table}[h] 
\footnotesize
\centering
\setlength{\extrarowheight}{2pt}
	\begin{tabular}{| l | l | >{\centering\arraybackslash}p{1.7cm}  | >{\centering\arraybackslash}p{1.7cm}  | >{\centering\arraybackslash}p{1.7cm}  | >{\centering\arraybackslash}p{1.7cm} |  } \hline 
		\multicolumn{2}{|c|}{$h$} & $1/20$  & $1/40$  & $1/80$  & $1/160$ \\ \hline
		\multirow{2}{*}{$L^{2}$ errors} &   $c$  &  0.0641 & 0.0321 	[1.00] &0.0160	[1.00] & 0.0080	[1.00]  \\ \cline{2-6}
		& $\bphi$ & 0.0453  & 0.0227 [1.00] &  0.0114 [0.99] & 0.0057 [1.00] \\ \hline
		\multirow{3}{*}{Jacobi} & 2-sided Robin & 21 & 21 & 23 & 25 \\ \cline{2-6}
		& 1-sided Ventcel & 11  & 11 & 12 & 13 \\  \cline{2-6}
		& weighted Ventcel & 11  & 11 & 12 & 13  \\ \hline
		\multirow{3}{*}{GMRES} & 2-sided Robin & 16 & 16 & 20 & 22 \\ \cline{2-6}
		& 1-sided Ventcel & 10 & 11  &11 &12 \\  \cline{2-6}
		& weighted Ventcel & 8 & 10 &  10 & 11 \\ \hline
	\end{tabular}  \vspace{0.2cm}
	\caption{[Test case 1] Accuracy in space (the convergence rates are shown in square brackets) and numbers of Jacobi and GMRES iterations for different optimized parameters.} 
\label{tab:test1accuSpace} \vspace{-0.2cm}
\end{table}
%
% Convergence in time
%
Next we fix $h=1/200$ and decrease the time step sizes with $\Delta t_{1} =\sfrac{3}{4}\Delta t_{2}$ to verify the order of accuracy in time.  In Table~\ref{tab:test1accuTime}, we show the relative $L^{2}(\Omega)-$norm errors of $c$ and $\bphi$ at $T=1$ by the Ventcel-OSWR method,  and numbers of iterations for different algorithms.  The results show that the accuracy in time is preserved with nonconforming time grids.  The convergence with Jacobi iteration is independent of the time step size while GMRES, though converges faster,  slightly depends on~$\Delta t$.  Again, using optimized Ventcel parameters reduces the number of iterations by nearly half (compared to using optimized Robin parameters).
%
% PURE DIRICHLET BCs
%
%\begin{table}[h] 
%\centering
%\setlength{\extrarowheight}{2pt}
%\begin{tabular}{| l | c | c | c | c | } \hline 
%		\multirow{2}{*}{$\Delta t_{2}$} & \multicolumn{2}{c|}{$L^{2}$ errors} & \multicolumn{2}{c|}{Number of iterations}\\\cline{2-3} \cline{4-5}
%		& $c$ & $\bphi$ & Jacobi & GMRES \\ \hline
%		%
%		$T/6$ &  0.1859	\phantom{[1.00]}&  0.2008 	\phantom{[1.00]}& 17  &  11\\
%		$T/12$ & 0.0708	[1.39] & 0.0768 [1.39]& 17 & 12\\
%		$T/24$ & 0.0301 	[1.23] & 0.0325 [1.24]& 17 & 13\\
%		$T/48$ & 0.0145 	[1.05] & 0.0150 [1.12]& 17 & 14\\ 
%		%$T/96$ & 0.0088 	[0.72] & 0.0080 [0.91]& 18 & 16\\ 
%		\hline
%	\end{tabular}  
%	\caption{[Test case 1] Accuracy in time, the convergence rates are shown in square brackets.} 
%\label{tab:test1accuTime}
%\end{table}
% a=[0.1859 0.0708 0.0301 0.0145 0.0088];
% b=[0.2008 0.0768 0.0325 0.0150 0.0080];

\begin{table}[h] 
\footnotesize
\centering
\setlength{\extrarowheight}{2pt}
	\begin{tabular}{| l | l | >{\centering\arraybackslash}p{1.7cm}  | >{\centering\arraybackslash}p{1.7cm}  | >{\centering\arraybackslash}p{1.7cm}  | >{\centering\arraybackslash}p{1.7cm} | } \hline 
		\multicolumn{2}{|c|}{$\Delta t_{2}$} & $T/6$  & $T/12$  & $T/24$  & $T/48$ \\ \hline
		\multirow{2}{*}{$L^{2}$ errors} &   $c$  &  0.1859  & 0.0708	[1.39] &  0.0301 	[1.23]  &  0.0145 	[1.05]   \\ \cline{2-6}
		& $\bphi$ &  0.2008 & 0.0768 [1.39] & 0.0325 [1.24] & 0.0150 [1.12] \\ \hline
		\multirow{3}{*}{Jacobi} & 2-sided Robin & 33  & 33 & 33 & 35  \\ \cline{2-6}
		& 1-sided Ventcel & 17  & 17 & 17 & 17 \\  \cline{2-6}
		& weighted Ventcel & 17 & 17 & 17  & 17   \\ \hline
		\multirow{3}{*}{GMRES} & 2-sided Robin & 18 & 18 & 20 & 24 \\ \cline{2-6}
		& 1-sided Ventcel & 11 & 12  &13 &14 \\  \cline{2-6}
		& weighted Ventcel & 10  & 11 & 12 & 13 \\ \hline
	\end{tabular}  \vspace{0.2cm}
	\caption{[Test case 1] Accuracy in time (the convergence rates are shown in square brackets) and numbers of Jacobi and GMRES iterations for different optimized parameters.} 
\label{tab:test1accuTime} \vspace{-0.2cm}
\end{table}

\subsection{Test case 2: with piecewise discontinuous coefficients}
We now analyze the convergence of the iterative algorithms when the physical coefficients are discontinuous across the interface.  Towards that end, we consider the error equation with the same two nonoverlapping subdomains as in Test case~1. The porosity is $ \omega _{1} = \omega_{2} =\omega = 1 $. The diffusion and advection coefficients, $ d_{i} $ and $ \bu_{i} $,  for $ i~=1,2, $ are given in Table~\ref{tab:Test2Discont} for the diffusion-dominant, mixed regime and advection-dominant problems, respectively.  Note that the global P\'eclet number in each subdomain is computed by $$ \text{Pe}_{G,i}:=\frac{H\mid \bu_{i} \mid}{d_{i}}, \; i=1,2,$$ where H is the size of the domain (in this case, $ H=1 $). 

\begin{figure}[!http]
\centering
    \begin{subfigure}[b]{1 \textwidth}
    \centering
        \includegraphics[width=0.48\textwidth]{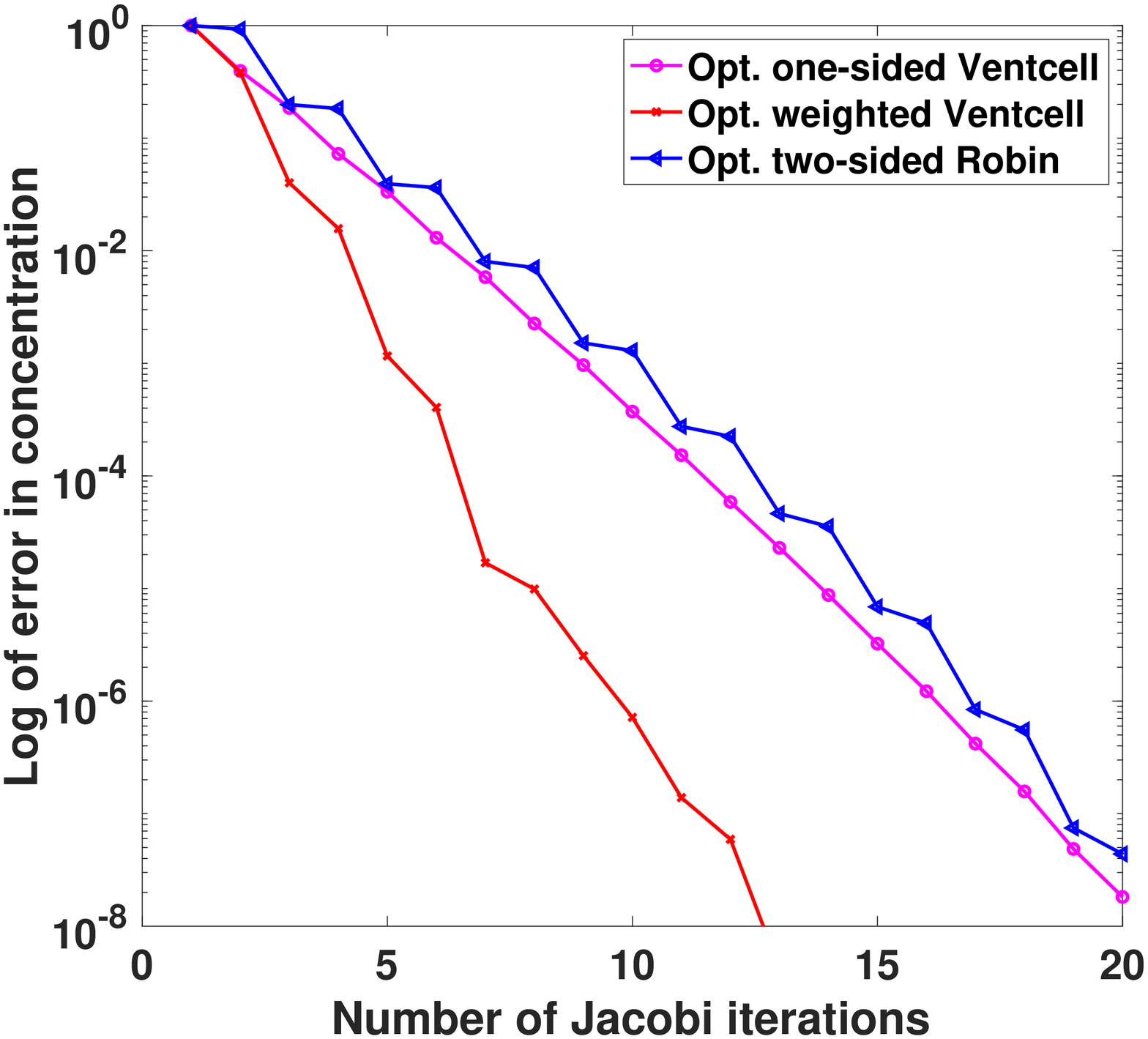}
         \includegraphics[width=0.48\textwidth]{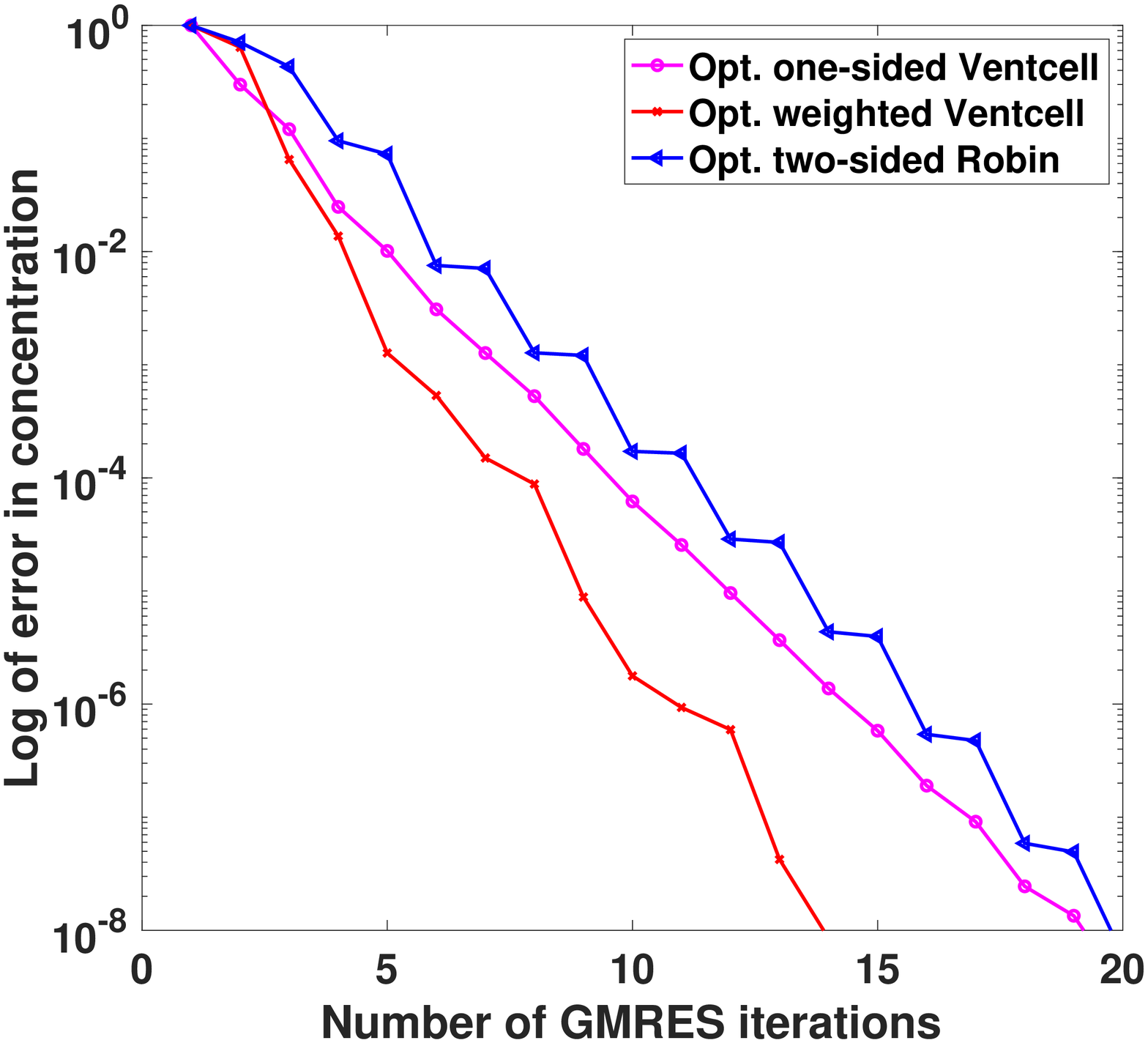}
        \caption{Diffusion dominance,  $\max Pe_{G,i} \approx 0.5$}
    \end{subfigure}\hspace{-0.25cm}
    \begin{subfigure}[b]{1\textwidth}
     \centering
        \includegraphics[width=0.48\textwidth]{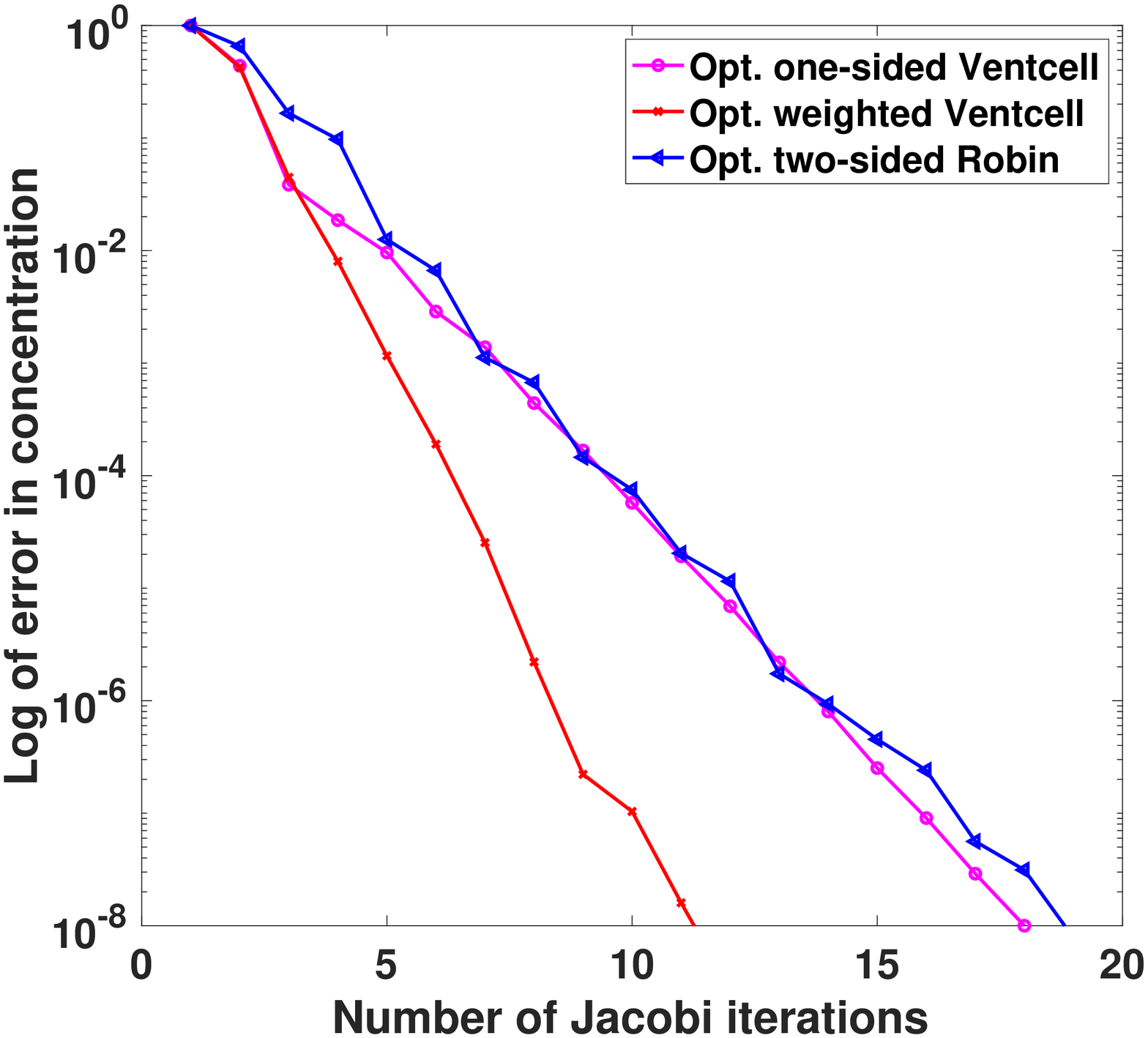}
         \includegraphics[width=0.48\textwidth]{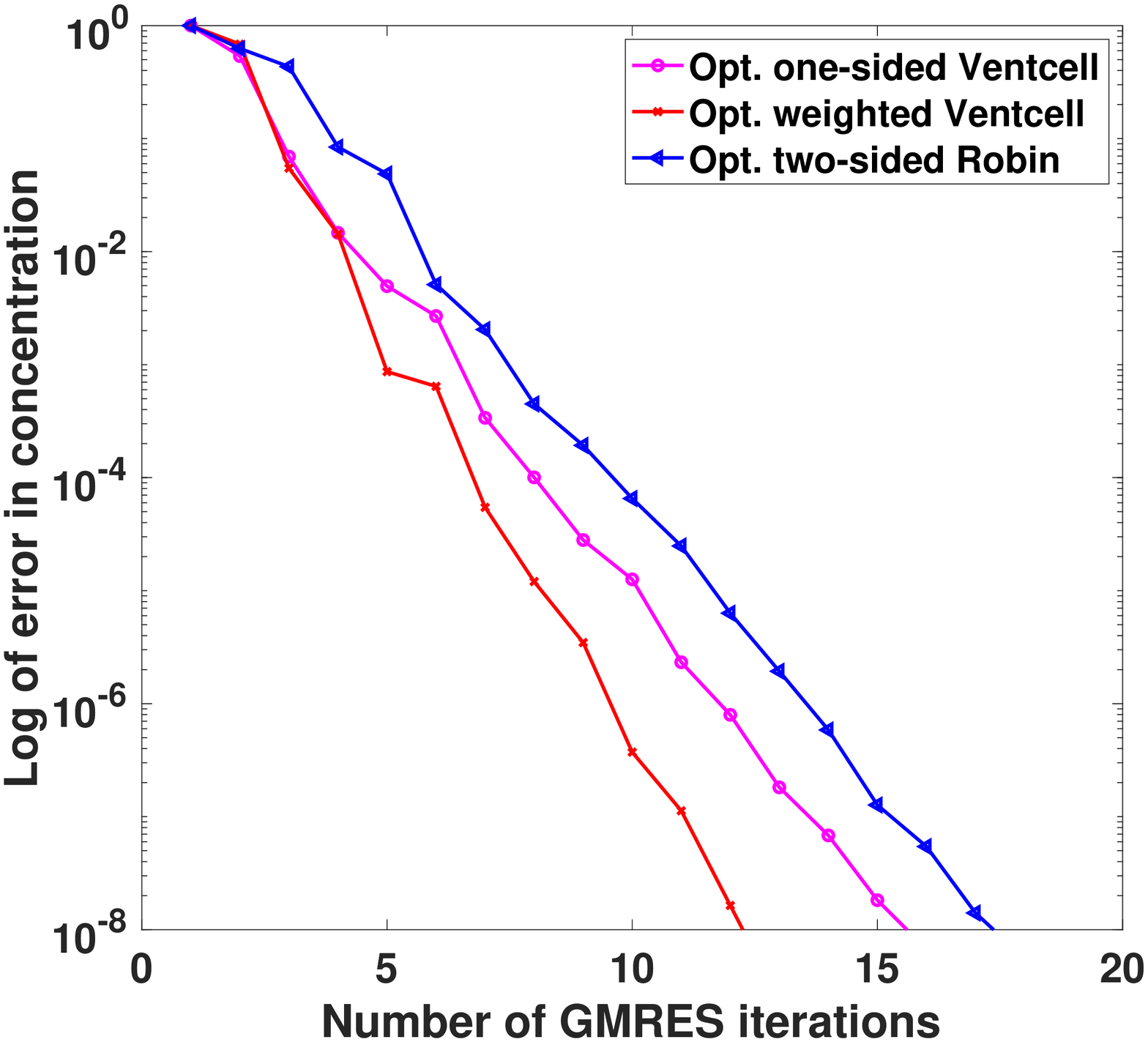}
        \caption{Mixed regime,  $\max Pe_{G,i} \approx 50$}
    \end{subfigure}\hspace{-0.25cm}
    \begin{subfigure}[b]{1\textwidth}
     \centering
        \includegraphics[width=0.48\textwidth]{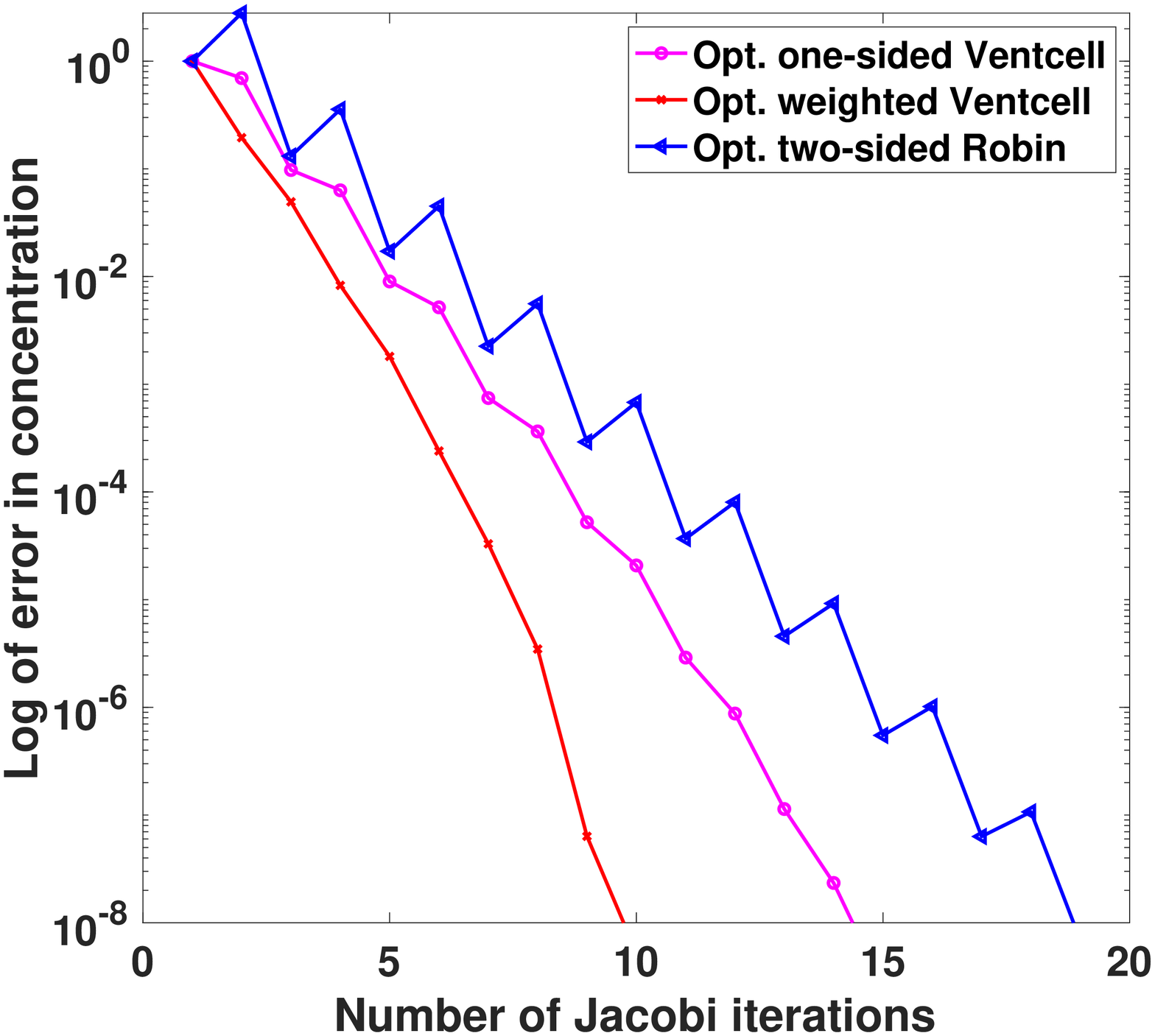}
         \includegraphics[width=0.48\textwidth]{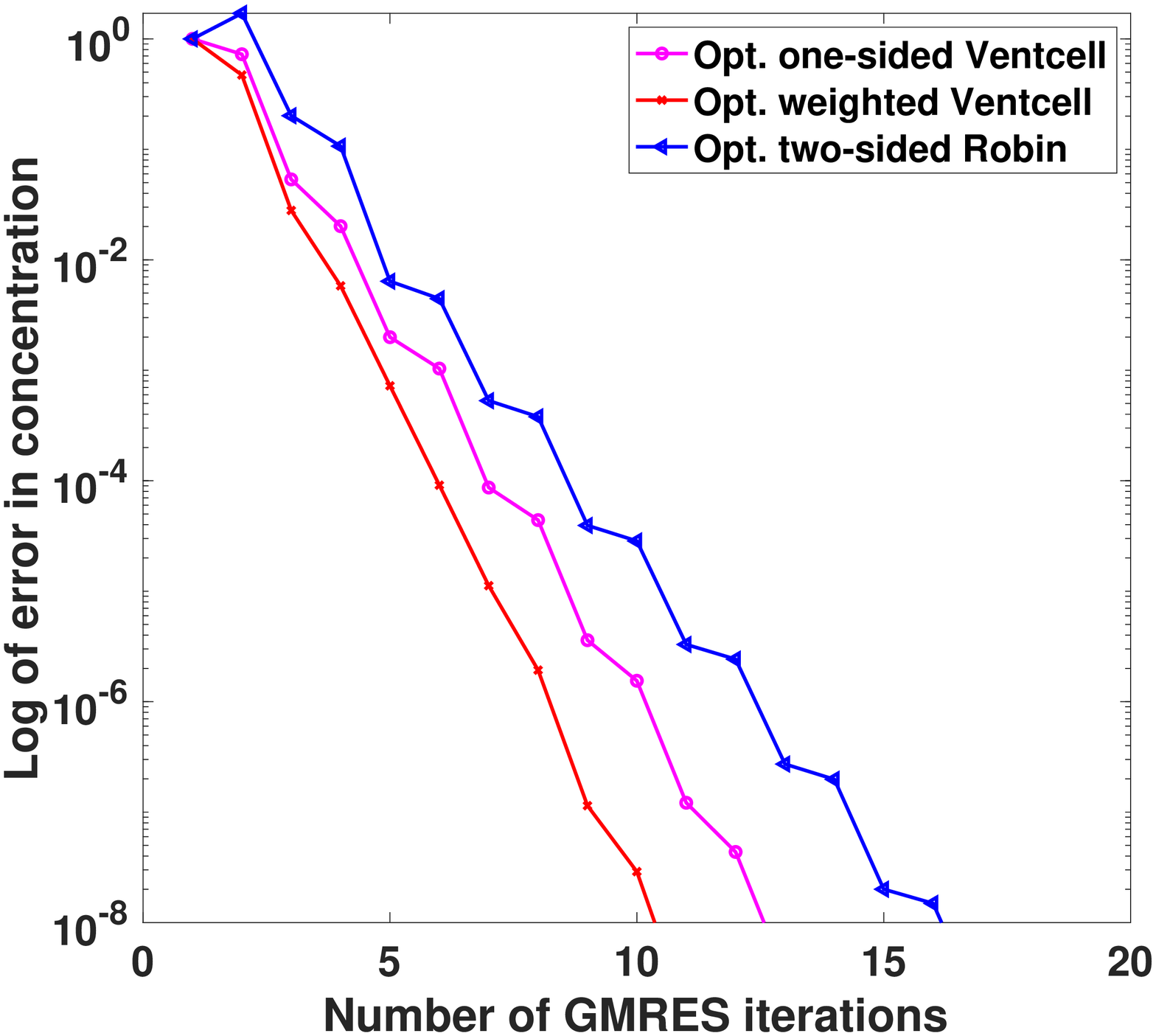}
        \caption{Advection dominance,  $\max Pe_{G,i} \approx 255$}
    \end{subfigure} \vspace{-0.2cm}
    \caption{[Test case 2] Convergence curves by Jacobi (left) and GMRES (right) for different P\'eclet numbers: $L^{2}-$norm errors in the concentration at $T=1$ with optimized two-sided Robin (blue curves),  optimized one-sided Ventcel (magenta curves), and optimized weighted Ventcel (red curves) parameters.}\label{fig:test2GMRESconv} 
\end{figure}

%In space, we use a conforming rectangular mesh $h= 1/100 $; in time, nonconforming time grids are considered with $\Delta t_{1}~=\sfrac{3}{4} \Delta t_{2}$ and $\Delta t_{2}=1/75$.
\begin{table}[h]
\footnotesize
\centering
\begin{tabular}{|l|l|l|l|l|l|l|}
  \hline
  		Problems							&$  d_{1}$	 		& $ \bu_{1} $		& $ \text{Pe}_{G,1}$ &$  d_{2}$	 & $ \bu_{2} $				& $ \text{Pe}_{G,2}$ 	 \\ \hline
   (a) Diffusion dominance    &  $1$  &   	$(-0.02, \; -0.5)^{T}$ &  $\approx 0.5$  & $0.1$ & $(-0.02, \; -0.05)^{T}$ & $\approx 0.5$	\\ \hline % 2c
  (b) Mixed regime     &  $0.01$  &   	$(-0.02, \; -0.5)^{T}$ &  $\approx 50$  & $0.1$ & $(-0.02, \; -0.05)^{T}$ & $\approx 0.5$		\\ \hline % 2a
  (c) Advection dominance     &  $0.02$  &   	$(0.5, \; 1)^{T}$ &  $\approx 56$  & $0.002$ & $(0.5, \; 0.1)^{T}$ & $\approx 255$		\\ \hline % 2b
\end{tabular} \vspace{0.2cm}
\caption{[Test case 2] Discontinuous diffusion and advection coefficients.}
\label{tab:Test2Discont} \vspace{-0.2cm}
\end{table}
We consider $h=1/100$,  $\Delta t_{1} =1/100$ and $\Delta t_{2}=1/75$, and use a random initial guess on the space-time interface to start the iteration.  Figure~\ref{fig:test2GMRESconv} show the errors (in logarithmic scale) in $L^{2}(\Omega)-$norm of the concentration versus the number of Jacobi or GMRES iterations (similar convergence curves are obtained for the vector variable and are omitted).  Three choices of optimized parameters are considered: optimized two-sided Robin (blue curves),  optimized one-sided Ventcel (magenta curves), and optimized weighted Ventcel (red curves).  We observe that for coefficients with jumps,  optimized weighted Ventcel parameters are robust with respect to different Pecl\'et numbers and give faster convergence than optimized one-sided Ventcel parameters. The one-sided parameters has similar performance as the two-sided Robin parameters when advection is not so strong.  We see that GMRES slightly improves the convergence when compared to Jacobi; however,  there is not much difference between GMRES and Jacobi when optimized weighted Ventcel parameters. 

\begin{figure}[!http]
\centering
    \begin{subfigure}[b]{0.48\textwidth}
    \centering
        \includegraphics[width=\textwidth]{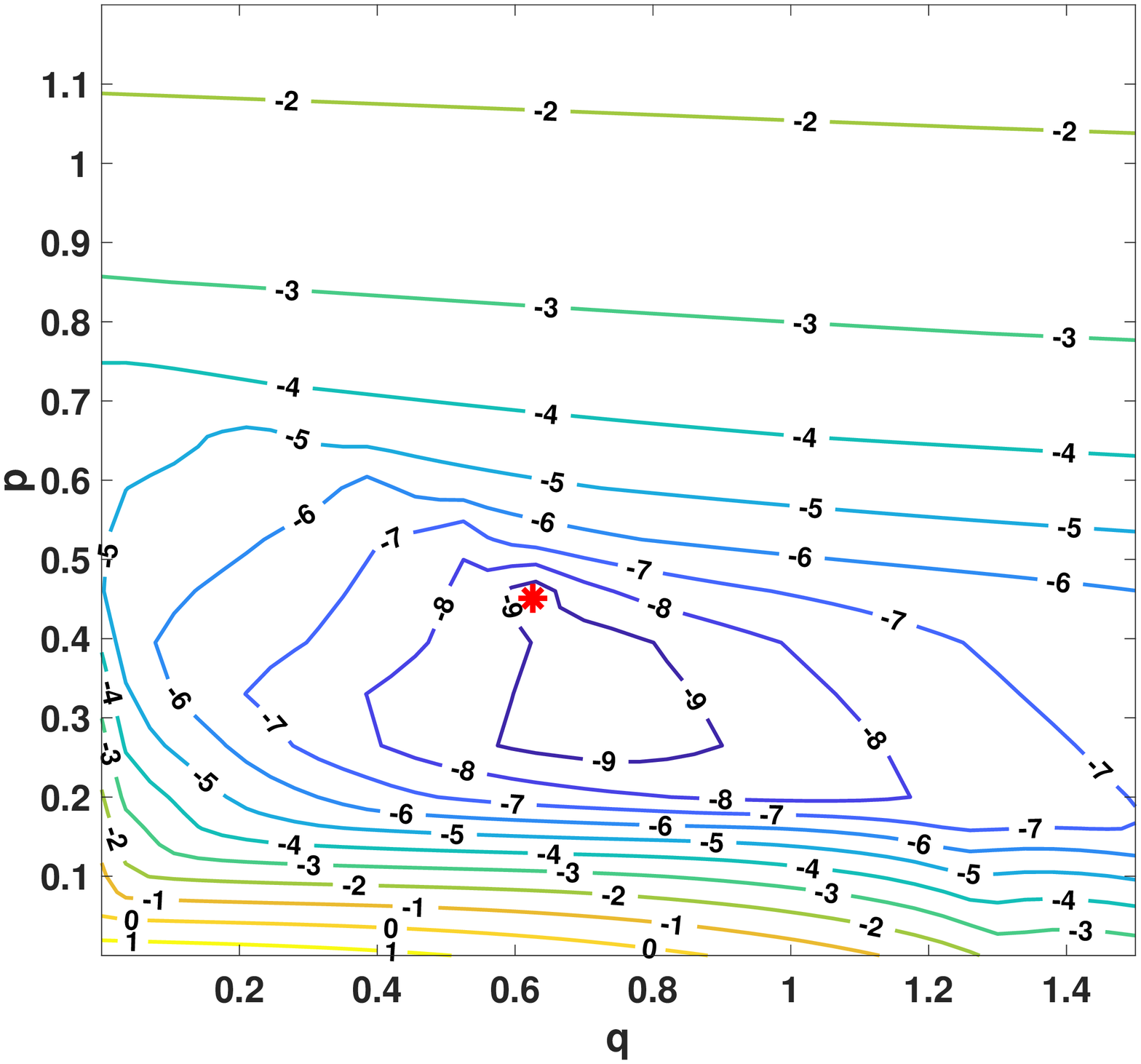}
    \end{subfigure} \hspace{4pt}
    \begin{subfigure}[b]{0.48\textwidth}
     \centering
        \includegraphics[width=\textwidth]{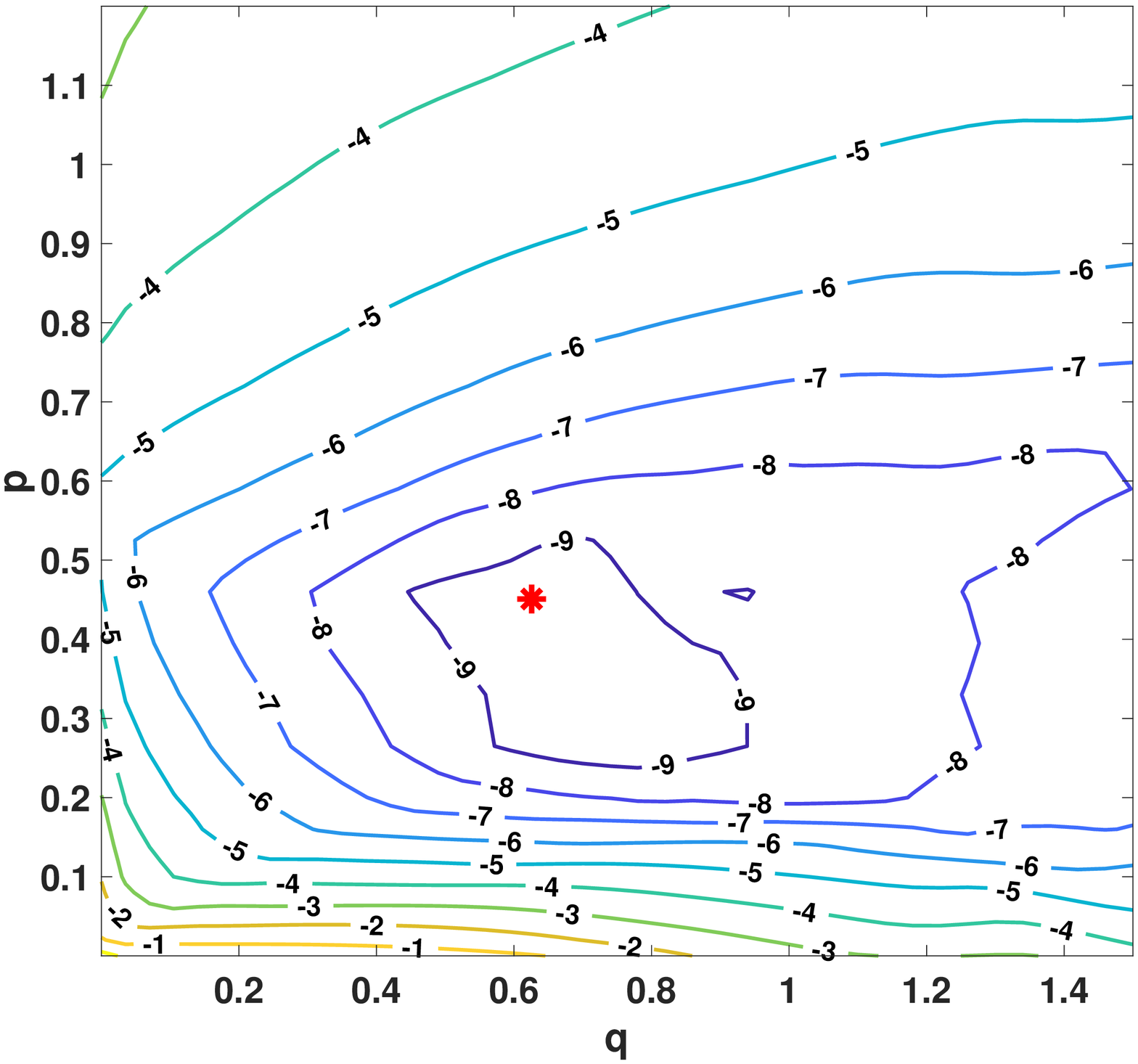}
    \end{subfigure} \vspace{-0.2cm}
    \caption{[Test case 2: Advection dominance] Level curves for the error in concentration (in logarithmic scales) after 12 iterations of Jacobi (left) and GMRES (right) for various values of $p$ and $q$. The red star shows the optimized values computed by numerically minimizing the continuous convergence factor of the OSWR algorithm. }\label{fig:test2Voptpara}
\end{figure}

In the following, we shall use only optimized weighted Ventcel parameters since they give fastest convergence.  To verify the performance of the optimized parameters,  we consider the advection-dominant problem (i.e. Problem (c)) and show in Figure~\ref{fig:test2Voptpara} the errors in concentration (in logarithmic scale) for various values of the parameters $p$ and $q$ after 12 Jacobi and GMRES iterations.  We see that for both Jacobi and GMRES,  the pair of optimized parameters $(q^{\ast},p^{\ast})$ (red star) is located close to those giving the smallest error after the same number of iterations.

\begin{figure}[h]
\centering
    \begin{subfigure}[b]{0.48\textwidth}
    \centering
        \includegraphics[width=\textwidth]{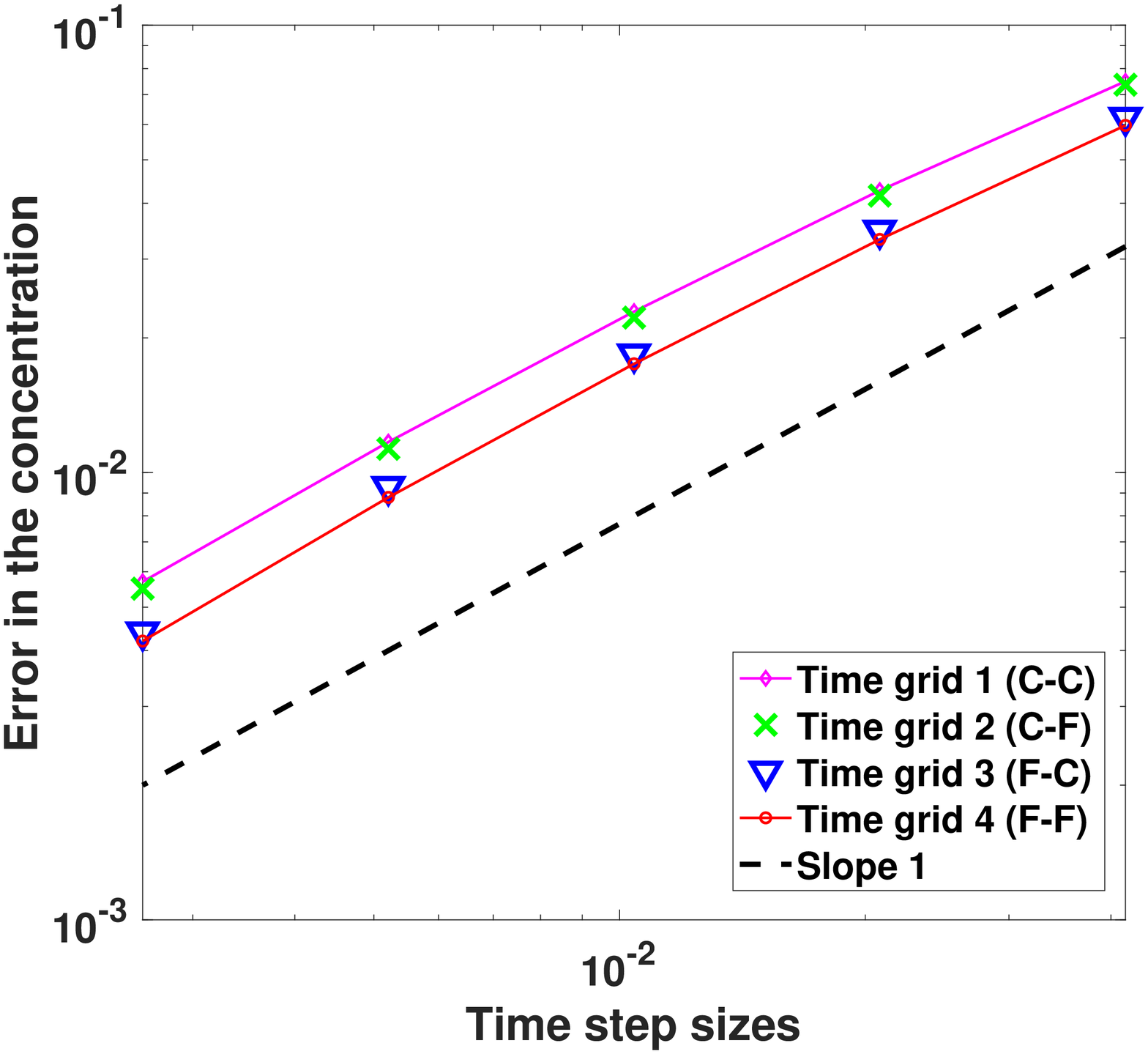}
    \end{subfigure} \hspace{4pt}
    \begin{subfigure}[b]{0.48\textwidth}
     \centering
        \includegraphics[width=\textwidth]{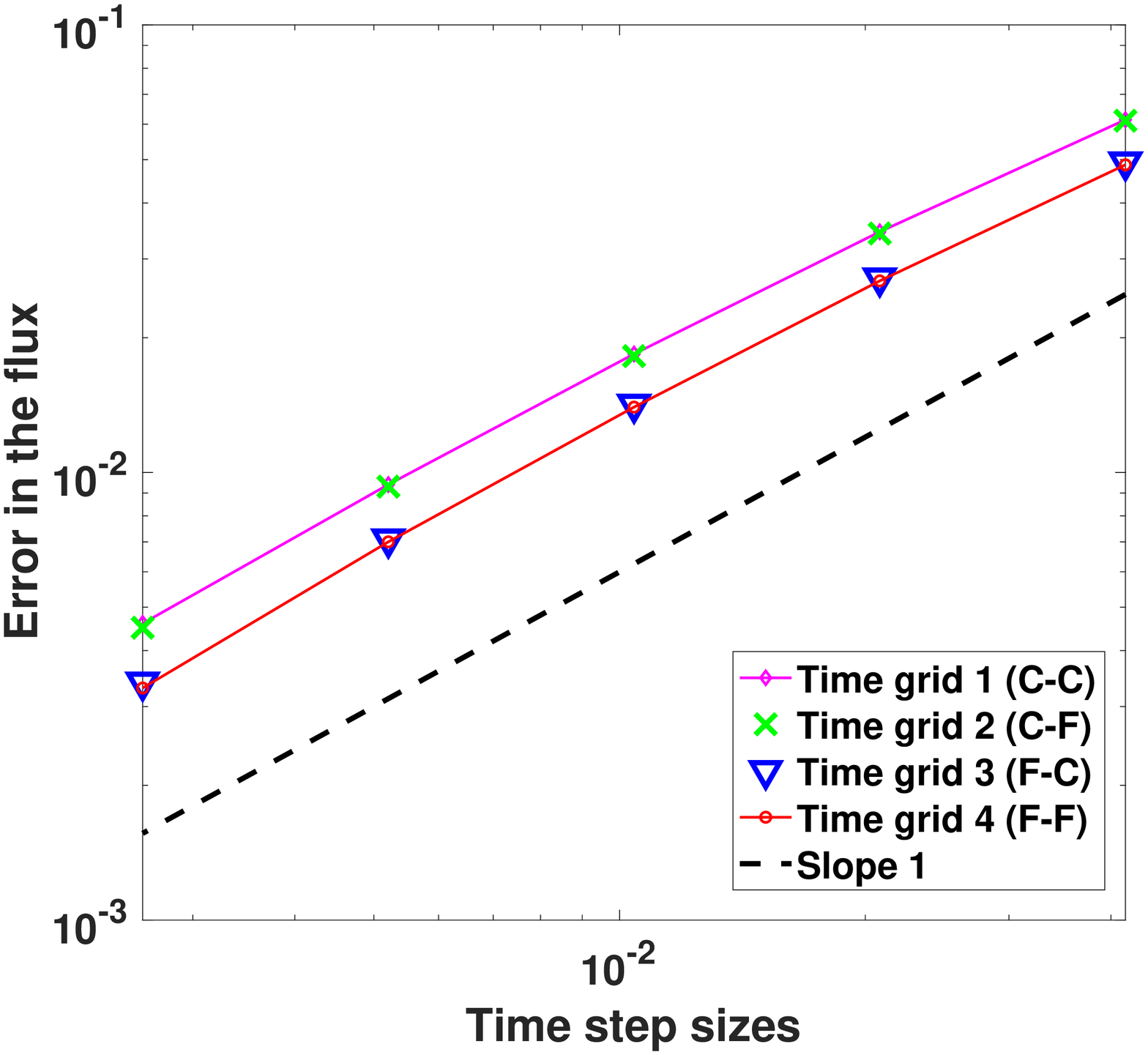}
    \end{subfigure} \vspace{-0.2cm}
    \caption{[Test case 2] Errors in the concentration $c$ (left) and the vector field $\bphi$ (right) between the reference and multidomain solutions.}\label{fig:test2accuracytime} 
\end{figure}

Finally,  we investigate whether the nonconforming time grids preserve the accuracy in time.  We again consider the  advection-dominant problem with homogeneous Dirichlet conditions on the boundary.  The source term is $ f(x,y, t)=\exp(-100((x-0.2)^2+(y-0.2)^2)),$ and the initial condition $p_{0}(x,y)=xy(1-x)(1-y)\exp(-100((x-0.2)^2+(y-0.2)^2)).$ We use four initial time grids with $ \Delta t_{c}=T/12$ and $ \Delta t_{f}=T/16 $ where $T=0.5$:
\begin{itemize} \itemsep0pt
	\item Time grid 1 (coarse-coarse): conforming with $ \Delta t_{1} = \Delta t_{2} = \Delta t_{c} $.
	\item Time grid 2 (coarse-fine): nonconforming with $ \Delta t_{1} = \Delta t_{c} $ and
              $ \Delta t_{2} = \Delta t_{f} $.
	\item Time grid 3 (fine-coarse): nonconforming with $ \Delta t_{1} = \Delta t_{f} $ and
              $ \Delta t_{2} = \Delta t_{c} $.
	\item Time grid 4 (fine-fine): conforming with $ \Delta t_{1} = \Delta t_{2} = \Delta t_{f} $.
\end{itemize}
The time steps are then refined several times by a factor of 2. In space, we fix a conforming rectangular mesh with $h=1/200$, and we compute a reference solution by solving the monodomain problem~\eqref{eq:mixedweak} directly
on a very fine time grid, with $ \Delta t = \Delta t_{f}/ 2^{7} $. The converged DD solution is such that the relative residual is smaller than $ 10^{-8} $.  We show in Figure~\ref{fig:test2accuracytime} the relative errors at $T=0.5$ versus the time step $ \Delta t = \max (\Delta t_{c}, \Delta t_{f}) $ obtained by using the optimized two-sided Ventcel conditions.  
We observe that first order convergence is preserved in the nonconforming case.  The errors obtained in the
nonconforming case with a fine time step in $\Omega_{1}$ where the parameters are large (Time grid 3 with blue triangle markers) are nearly the same as in the finer
conforming case (Time grid 4, in red with circle markers).  
On the other hand,  the errors obtained in the
nonconforming case with a fine time step in $\Omega_{2}$ where the parameters are small (Time grid 2 with green x-markers) are close to those by the coarse
conforming case (Time grid 1, in magenta with diamond markers).  Thus using nonconforming grids can
adapt the time steps in the subdomains depending on the physical parameters and limit the computational cost locally,  while preserving almost the same accuracy as in the finer conforming case.

\section*{Conclusion} We have developed mixed formulations of a global-in-time domain decomposition method based on OSWR and Ventcel transmission conditions to solve heterogeneous, linear advection-diffusion problems. A flux variable representing the total advective and diffusive flux is introduced,  and Lagrange multipliers are considered on the interfaces of the nonoverlapping subdomains to handle tangential derivatives involved in Ventcel conditions.  A space-time interface problem is derived based on trace operators of Ventcel-to-Ventcel type which involve solving time-dependent problems with Ventcel boundary conditions in the subdomains.  The local Ventcel problem is written in mixed form and can be seen as a coupling of a $d$-dimensional PDE in the subdomain  and a $(d-1)$-dimensional PDE on the interface.  Mixed hybrid finite elements are used to discretize the equations in space, where the Lagrange multiplier arising in the hybridization is used to discretize the advective term.  The fully discrete problems are obtained by the backward Euler method. The proposed global-in-time optimized Ventcel-Schwarz method is fully implicit and allows different time steps in the subdomains.  Numerical results show that the method preserves the accuracy in time with nonconforming time grids; its convergence with optimized Ventcel transmission conditions is almost independent of the mesh size and time step size,  and is robust for problems with high Pecl\'et numbers.  Moreover,  using optimized weighted Ventcel parameters for problems with discontinuous coefficients is most effective,  with a reduction of the number of iterations by nearly half compared to using optimized Robin parameters.  Future work includes the analysis of the discrete local Ventcel problems and convergence of the iterative algorithm, as well as the handling of Ventcel transmission conditions at cross points.

\appendix

\section{Algebraic system for the local Ventcel problem} \label{app:algsystem} We derive the linear system for the mixed Ventcel subdomain problem presented in Subsection~\ref{subsec:semilocV}; such a system is local to each subdomain and is solved independently at each Jacobi or GMRES iteration of the interface problem.  In the following,  the subscript $i$ refers to the $i$th subdomain.  We assume the porosity is piecewise constant and denote by $\omega_{i,K}$ the value of $\omega_{i}$ on element $K \in \iK_{h,i}$.  Taking $\mu=\chi_{K}$ ($K \in \iK_{h,i}$) in \eqref{eq:semiVmass},  $\bv=\bw_{KE}$ ($K \in \iK_{h,i}$, $E \subset \partial K$) in \eqref{eq:semiVflux} and $\eta=\chi_{E}$ ($E \in \iE_{h,i}^{I}$) in \eqref{eq:semiVlag} we have:
\begin{align}
&\vert K \vert \omega_{i,K} \partial_{t} c_{i,K} + \sum_{E \subset \partial K} \varphi_{i,KE} =\int_{K} f \, d\bx, \quad \forall K \in \iK_{h,i},  \label{eq:disVmass}\vspace{3pt}\\
&\begin{array}{c}
\displaystyle
\sum_{E^{\prime} \subset \partial K} A_{KEE^{\prime}} \varphi_{i,KE^{\prime}} - \sum_{E^{\prime} \subset \partial K} A_{KEE^{\prime}} u_{i,KE^{\prime}} \; \theta_{i,E^{\prime}}  -c_{i,K} + \theta_{i,E}=0, \\
\hspace{5cm}\forall K \in \iK_{h,i}, \forall E \subset \partial K, 
\end{array} \label{eq:disVflux}\vspace{3pt}\\
& \varphi_{i,K_{1}E} +  \varphi_{i,K_{2}E}  = 0, \quad \forall E=\left (\partial K_{1} \cap \partial K_{2}\right ) \in \iE_{h,i}^{I},  \; K_{1}, K_{2} \in \iK_{h,i},  \label{eq:disVlag}
\end{align}
where $A_{KEE^{\prime}} = \int_{K} (\bD^{-1} \bw_{KE^{\prime}}) \cdot \bw_{KE} \, d\bx$, for $E, E^{\prime} \subset \partial K$, $K \in \iK_{h,i}$. 

For each interface edge $E \in \iE_{h}^{\Gamma}$,  we write $E=[P_{1}, P_{2}]$ where $P_{1}, P_{2} \in \iP_{h}^{\Gamma}$ are the two endpoints of $E$.  For the equations on the interface,  we proceed in a similar manner, i.e.  choosing the test functions in \eqref{eq:semiV1dmass}, \eqref{eq:semiV1dflux} and \eqref{eq:semiV1dlag} to be the basis functions of the associated spaces, and obtain the following equations after some simple calculations:
\begin{align}
&\begin{array}{c}
\displaystyle
-\varphi_{i,KE} + \alpha_{i,j} \lambda_{i,E} + \beta_{i,j} \omega_{i,K} \partial_{t} \lambda_{i,E} + \beta_{i,j} \left (\varphi_{\Gamma,i,EP_{1}}+\varphi_{\Gamma,i,EP_{2}}\right ) = \int_{E} \zeta dy,  \vspace{4pt}\\
\hspace{5cm} \forall E=[P1,P2] \in \iE_{h}^{\Gamma},  E\subset \partial K \in \iK_{h,i}, 
\end{array} \label{eq:disV1dmass} \vspace{3pt}\\
& \frac{1}{6}\bD_{j,\Gamma}^{-1}\vert E \vert \left [ \begin{array}{cc} 2 & -1  \\ -1  & 2 \end{array} \right ]\, \left [ \begin{array}{c} \varphi_{\Gamma,i,EP_{1}}  \\ \varphi_{\Gamma,i,EP_{2}} \end{array} \right ] - \frac{1}{6}\bD_{j,\Gamma}^{-1}\vert E \vert \left [ \begin{array}{cc} 2 & -1  \\ -1  & 2 \end{array} \right ] \, \left [ \begin{array}{c} u_{i,\Gamma,EP_{1}} \; \xi_{\Gamma,i,P_{1}}  \\ u_{i,\Gamma,EP_{2}} \; \xi_{\Gamma,i,P_{2}} \end{array} \right ] \nonumber \vspace{6pt}\\
&\hspace{1.5cm}-\lambda_{i,E} \left [ \begin{array}{c} 1 \\ 1 \end{array} \right ] + \left [ \begin{array}{c} \xi_{\Gamma,i,P_{1}} \\ \xi_{\Gamma,i,P_{2}} \end{array} \right ] =\left [ \begin{array}{c} 0 \\ 0 \end{array} \right ], \; \forall E=[P_{1},P_{2}] \in \iE_{h}^{\Gamma},  \label{eq:disV1dflux}\vspace{5pt}\\ 
&\varphi_{\Gamma,i,E_{1}P}+\varphi_{\Gamma,i,E_{2}P}=0, \quad \forall P = \partial E_{1} \cap \partial E_{2} \in \iP_{h}^{\Gamma}, E_{1}, E_{2} \in \iE_{h}^{\Gamma}.  \label{eq:disV1dlag}
\end{align}
Note that $\bD_{j,\Gamma}$ is a scalar for $\Gamma \subset \mR^{2}$.  With backward Euler time stepping (cf. Section~\eqref{sec:time}),  we obtain,  at each time step~$t_{i,m}$, for $m=1, \hdots, M_{i}$,  a symmetric system of the following form  : 
\begin{equation} \label{eq:linearsys}
\left [ \begin{array}{cccccc} A & B^{T} & E^{T} & E_{\Gamma}^{T} &  0 & 0  \\
B & C & 0 & 0 & 0 & 0 \\
E &0  & 0 & 0 & 0 & 0 \\
E_{\Gamma} & 0 & 0 & C_{\alpha,\beta} & B_{\Gamma} & 0 \\
0 & 0 & 0 & B_{\Gamma}^{T} & A_{\Gamma} & E_{\partial \Gamma}^{T} \\
0  & 0 & 0 & 0 & E_{\partial \Gamma} & 0 
\end{array} \right ] 
\left [ \begin{array}{c} \widetilde{\varphi_{i}}^{m} \\ \widetilde{c}_{i}^{m}\\ \widetilde{\lambda}_{i,I}^{m} \\ \widetilde{\lambda}_{i,\Gamma}^{m} \\ \widetilde{\varphi}_{\Gamma,i}^{m} \\ \widetilde{\xi}_{\Gamma,i}^{m} \end{array} \right ]  =  \left [ \begin{array}{c} 0 \\ \ast \\ 0 \\ \ast \\ 0 \\ 0 \end{array} \right ],
\end{equation}
where $\widetilde{\varphi}_{i}^{m}$ and $\widetilde{c}_{i}^{m}$ represent the flux and concentration 
unknowns in subdomain $\Omega_{i}$ at time $t_{i,m}$ (i.e.  $\left (\varphi_{i,KE}(t^{m}\right )_{K \in \iK_{h,i}, E \subset \partial K}$ and $\left (c_{i,K}(t^{m})\right )_{K\in \iK_{h,i}}$).  $\widetilde{\lambda}_{i,I}^{m}$ and $\widetilde{\lambda}_{i,\Gamma}^{m}$ are the Lagrange multipliers living on internal edges and on interface edges,  respectively.  $\widetilde{\varphi}_{\Gamma,i}^{m}$ represents the one-dimensional flux unknowns on interface and $\widetilde{\xi}_{\Gamma,i}^{m}$ is the Lagrange multipliers associated with the interface nodes.  The first two equations of~\eqref{eq:linearsys} correspond to the flux equation~\eqref{eq:disVflux} and the mass conservation equation~\eqref{eq:disVmass} while the third equation represents~\eqref{eq:disVlag} to enforce the continuity of the flux over internal edges. The last three equations of \eqref{eq:linearsys} correspond to \eqref{eq:disV1dmass}, \eqref{eq:disV1dflux} and \eqref{eq:disV1dlag} to enforce the Ventcel condition on the interface.  The right-hand side of \eqref{eq:linearsys} whose explicit form is omitted can be deduced from \eqref{eq:disVmass} and \eqref{eq:disV1dmass}. 

\section{Computing the discrete Ventcel-to-Ventcel operators} \label{app:comphi}
The interface problem~\eqref{eq:semiIP} involves the calculation of the Ventcel-to-Ventcel operators $\iS_{h,i}^{\text{VtV}}(\zeta,f,c_{0})$ after solving the time-dependent subdomain problem~\eqref{eq:disVmass}-\eqref{eq:disV1dlag} over the whole time interval $(0,T)$.  The interface space $\Theta_{h,\Gamma}$ consists of piecewise constant functions on the interface edges $E \in \iE_{h}^{\Gamma}$,  thus we choose $\theta_{i}=\chi_{E}$ and rewrite~\eqref{eq:semiIP} as
\begin{equation}
\int_{E} \zeta_{i} \, dy - \int_{E} \iS_{h,j}^{\text{VtV}}(\zeta_{j},0,0) \, dy=\int_{E} \iS_{h,j}^{\text{VtV}}(0,f,c_{0}) \, dy, \; \forall E \in \iE_{h}^{\Gamma}, 
\end{equation}
for a.e. $t \in (0,T)$ and for  $i=1,2, \, j=(3-i)$. Using \eqref{eq:semiSvtv}, we deduce that:
\begin{equation} \label{eq:disSvtv}
\begin{array}{ll}
\displaystyle
\int_{E} \iS_{h,i}^{\text{VtV}}(\zeta_{i},f,c_{0}) \, dy&=\bphi_{i,KE}+ \alpha_{j,i} \; \lambda_{i,E} + \beta_{j,i}  \omega_{i} \partial_{t}\lambda_{i,E} \\
&\hspace{1.5cm}+\beta_{j,i}\left (\phi_{i,\Gamma, EP_{1}}+\phi_{i,\Gamma, EP_{2}}\right), 
\end{array}
\end{equation}
for $E=[P_{1},P_{2}] \in \iE_{h}^{\Gamma}$,  where we have used the fact that $\bpsi_{i,\Gamma}(t) \in \Sigma_{h,\Gamma}$ with 
$$\bpsi_{i,\Gamma}(t,y)=\sum_{P \in \partial E} \phi_{i,\Gamma,EP} (t)\; \bw_{\Gamma,EP}(y).$$
We will compute the terms involving $\left (\phi_{i,\Gamma,EP}\right )_{E \in \iE_{h}^{\Gamma}, P \in \partial E}$ in \eqref{eq:disSvtv} by using the relation \eqref{eq:tantracedis} between $\bpsi_{i,\Gamma}$ and the solution of the subdomain problem~\eqref{eq:disVmass}-\eqref{eq:disV1dlag}.  Taking $\bv_{\Gamma}$ to be the basis functions of $\Sigma_{h,\Gamma}$ in  \eqref{eq:tantracedis} and after some simple calculations, we obtain
\begin{equation} \label{eq:tantracemat}
\begin{array}{ll}
\left [ \begin{array}{cc} 2 & -1  \\ -1  & 2 \end{array} \right ]\, \left [ \begin{array}{c} \phi_{i,\Gamma, EP_{1}}  \\ \phi_{i,\Gamma,EP_{2}} \end{array} \right ] &= \left [ \begin{array}{cc} 2 & -1  \\ -1  & 2 \end{array} \right ] \, \left [ \begin{array}{c} u_{i,\Gamma,EP_{1}} \; \xi_{\Gamma,i,P_{1}}  \\ u_{i,\Gamma,EP_{2}} \; \xi_{\Gamma,i,P_{2}} \end{array} \right ] \vspace{6pt}\\
&\hspace{-2.5cm}-\bD_{j,\Gamma}^{-1}\bD_{i,\Gamma}\left [ \begin{array}{cc} 2 & -1  \\ -1  & 2 \end{array} \right ] \, \left [ \begin{array}{c} u_{j,\Gamma,EP_{1}} \; \xi_{\Gamma,i,P_{1}}  \\ u_{j,\Gamma,EP_{2}} \; \xi_{\Gamma,i,P_{2}} \end{array} \right ] \vspace{6pt}\\
&\hspace{-2.5cm}+ \bD_{j,\Gamma}^{-1}\bD_{i,\Gamma} \left [ \begin{array}{cc} 2 & -1  \\ -1  & 2 \end{array} \right ]\, \left [ \begin{array}{c} \varphi_{\Gamma, i,EP_{1}}  \\ \varphi_{\Gamma,i,EP_{2}} \end{array} \right ],  \; \forall E=[P_{1},P_{2}] \in \iE_{h}^{\Gamma},
\end{array}
\end{equation} 
Adding the two equations of \eqref{eq:tantracemat} yields the sum of $\phi_{i,\Gamma, EP_{1}}$ and $\phi_{i,\Gamma, EP_{2}}$ for each $E=[P_{1},P_{2}] \in \iE_{h}^{\Gamma}$. Then we plug it into \eqref{eq:disSvtv} to compute the Ventcel data for the interface problem.  Similar calculation can be done for the fully discrete case, in particular we have:
\begin{align*}
\int_{J_{i,m}}\int_{E} \iS_{h,i}^{\text{VtV}}(\zeta_{i},f,c_{0})\, dy&=\bphi_{i,KE}^{m} + \alpha_{j,i} \; \lambda_{i,E}^{m} + \beta_{j,i} \omega_{i} \left (\frac{\lambda_{i,E}^{m}-\lambda_{i,E}^{m-1}}{\Delta t_{i,m}} \right )\\
& \hspace{1.4cm} + \beta_{j,i} \left (\phi_{i,\Gamma, EP_{1}}^{m}+\phi_{i,\Gamma, EP_{2}}^{m}\right ),
\end{align*}
for $E=[P_{1},P_{2}] \in \iE_{h}^{\Gamma}$ and $m=1, \hdots, M_{i}$. 

\bibliographystyle{elsarticle-num}

\end{document}